\documentclass[a4paper,11pt]{amsart}
\usepackage{amsmath,amssymb,amsthm}
\usepackage{geometry}
\usepackage[dvipdfmx]{graphicx}
\allowdisplaybreaks

\newtheorem{thm}{Theorem}[section]
\newtheorem{lemma}[thm]{Lemma}
\newtheorem{theorem}[thm]{Theorem}
\newtheorem{proposition}[thm]{Proposition}
\newtheorem{corollary}[thm]{Corollary}

\newtheorem{theoremintro}{Theorem}
\newtheorem{corintro}[theoremintro]{Corollary}

\theoremstyle{definition}
\newtheorem{definition}[thm]{Definition}
\newtheorem{example}[thm]{Example}
\newtheorem{remark}[thm]{Remark}

\newcommand{\RR}{\mathbb R}
\newcommand{\RRR}{\mathbb {R}_{\geq 0}}

\newcommand{\NN}{\mathbb N}
\newcommand{\HH}{\mathbb H}
\renewcommand{\H}{\mathcal H}
\newcommand{\ZZ}{\mathbb Z}
\newcommand{\F}{\mathcal F}

\newcommand{\B}{\mathcal B}

\newcommand{\gL}{\Lambda}
\newcommand{\gG}{\Gamma}
\newcommand{\gS}{\Sigma}
\newcommand{\id}{\mathrm{id}}

\newcommand{\Sub}{\mathrm{Sub}}
\renewcommand{\:}{\colon}
\renewcommand{\tilde}{\widetilde}

\newcommand{\SC}{\mathrm{SC}}
\newcommand{\Stab}{\mathrm{Stab}}

\newcommand{\GC}{\mathrm{GC}}

\newcommand{\Isom}{\mathrm{Isom}^+}

\newcommand{\Map}{\mathrm{Map}}
\newcommand{\CH}{\mathrm{CH}}
\newcommand{\ML}{\mathrm{ML}}
\newcommand{\Thu}{\mathrm{Thu}}
\newcommand{\supp}{\mathrm{supp}}
\newcommand{\Area}{\mathrm{Area}}

\begin{document}
\title{Counting subgroups via Mirzakhani's curve counting}

\author[D. Sasaki]{Dounnu Sasaki}
\address{Department of mathematics, Faculty of Science, Josai University, Keyakidai 1-1, Sakado-shi, Saitama, 350-0295, Japan}
\email{dsasaki@josai.ac.jp}

\begin{abstract}
Given a hyperbolic surface $\Sigma$ of genus $g$ with $r$ cusps, Mirzakhani proved that the number of closed geodesics of length at most $L$ and of a given type is asymptotic to $cL^{6g-6+2r}$ for some $c>0$.
Since a closed geodesic corresponds to a conjugacy class of the fundamental group $\pi_1(\Sigma )$, we extend this to the counting problem of conjugacy classes of finitely generated subgroups of $\pi_1(\Sigma )$.
Using `half the sum of the lengths of the boundaries of the convex core of a subgroup' instead of the length of a closed geodesic, we prove that the number of such conjugacy classes is similarly asymptotic to $cL^{6g-6+2r}$ for some $c>0$.
As a special case, these conjugacy classes can be interpreted as subsurfaces of $\Sigma$ via their convex cores, and the result can be viewed as counting subsurfaces of a given type.
Furthermore, we see that the above length measurement for subgroups is `natural' within the framework of subset currents, which serve as a completion of weighted conjugacy classes of finitely generated subgroups of $\pi_1(\Sigma )$.
\end{abstract}

\maketitle

\section{Introduction}

Let $\gS$ be an orientable hyperbolic surface of genus $g$, possibly with $r$ cusps, having finite area and without boundary. Let $\Map (\gS):=\mathrm{Homeo}^+(\gS)/\mathrm{isotopy}$ denote the mapping class group of $\gS$.
We study the asymptotic growth of the number of $\Map(\gS)$ orbits of a conjugacy class of a finitely generated subgroup of $\pi_1(\gS)$, analogous to the results of Mirzakhani for closed geodesics on $\gS$.

Mirzakhani \cite{Mir08, Mir16} proved that for any weighted multicurve $\gamma_0$ and any finite-index subgroup $\Gamma$ of $\Map(\gS)$, there exists a positive constant $\mathfrak{c}^\Gamma_{g,r}(\gamma_0)$ such that
\[ \lim_{L\rightarrow \infty}\frac{\# \{ \gamma \in \Gamma (\gamma_0 )\mid \ell (\gamma )\leq L\}}{L^{6g-6+2r}}=\mathfrak{c}^\Gamma_{g,r}(\gamma_0) m_\Thu (\ell^{-1}([0,1])),\]
where $\ell$ is the length functional and $m_\Thu$ is the Thurston measure on the space $\ML (\gS)$ of measured laminations on $\gS$.
A weighted multicurve on $\gS$ is a formal finite sum of positive weighted closed geodesics on $\gS$, which is regarded as an element of the space $\GC(\gS)$ of geodesic currents on $\gS$. For any $[\phi]\in \Map (\gS)$ and a closed geodesic $\gamma$ of $\gS$, $[\phi](\gamma)$ represents the closed geodesic free homotopic to $\phi(\gamma)$, and this action extends linearly. Then, $\gG(\gamma_0 )$ is the set consisting of all $[\phi ](\gamma_0)$ for $[\phi] \in \gG$.

Erlandsson-Souto \cite{ES22} extended this theorem into a general version where the length functional $\ell$ can be replaced with any positive homogeneous continuous function
$F$ on the space $\GC_K(\gS)$ of geodesic currents supported by some compact subset $K$ of $\gS$, which includes all closed geodesics in $\gG (\gamma_0)$.
The function $F$ is homogeneous if $F(c\mu )=cF(\mu)$ for any $\mu \in \GC_K (\gS)$ and $c\in \RRR$, and $F$ is positive on $\GC_K(\gS)$ if $F(\mu )>0$ for any non-zero $\mu\in \GC_K (\gS)$.

Since an unoriented closed geodesic corresponds to a conjugacy class of a cyclic subgroup of $\pi_1(\gS)$, we are going to extend the above results to the counting problem of \emph{conjugacy classes of finitely generated subgroups of $\pi_1(\Sigma )$.}
Our first main result can be stated as follows:
\begin{theoremintro}[See Theorem \ref{thm:counting subgroup use length}]\label{thmintro:1}
Let $\gG$ be a finite-index subgroup of $\Map (\gS)$. Let $H$ be a finitely generated subgroup of $\pi_1(\gS)$ such that its limit set contains at least two points.
Then, there exists a constant $\mathfrak{c}^\Gamma_{g,r}(H)$ such that
\[ \lim_{L\rightarrow \infty}\frac{\# \{ [H'] \in \Gamma ([H] )\mid \ell_\SC ([H'] )\leq L\}}{L^{6g-6+2r}}=\mathfrak{c}^\Gamma_{g,r}(H) m_\Thu (\ell^{-1}([0,1])),\]
where $\ell_\SC$ is the generalized length functional for subgroups. Explicitly, $\ell_\SC ([H])$ is half the sum of the lengths of boundary components of the convex core of $H$.

This formula remains valid if the single conjugacy class $[H]$ is replaced by any formal weighted finite sum
\[ a_1[H_1]+\cdots +a_m[H_m]\]
where $a_1,\dots,a_m > 0$ and $H_1,\dots,H_m$ are finitely generated subgroups of $\pi_1(\gS)$, each with a limit set containing at least two points.
\end{theoremintro}

Note that $\Map(\gS)$ acts on the set of conjugacy classes of finitely generated subgroups of $\pi_1(\gS)$, and that $\gG([H])$ denotes the orbit of $[H]$ under the action of $\gG$.
The above condition that the limit set of $H$ contains at least two points is equivalent to saying that $H$ is not generated by a single element whose representative is peripheral. The constant $\mathfrak{c}^\Gamma_{g,r}(H)$ is positive if and only if $H$ is not a finite-index subgroup of $\pi_1(\gS)$.

\subsection{Geometric interpretation of Theorem \ref{thmintro:1}}

We denote by $\HH$ the hyperbolic plane. The fundamental group $\pi_1(\gS)$ acts on $\HH$ by deck transformations.
For any non-trivial finitely generated subgroup $H$ of $\pi_1(\gS)$ whose limit set contains at least two points, the convex core $C_H$ is defined as the smallest closed convex subset of the quotient space $\HH/H$ such that the inclusion map is a homotopy equivalence.
Let $p_H$ denote the restriction of the covering map $\HH/H\to \gS$ to $C_H$.

If a subgroup $H'$ of $\pi_1(\gS)$ is conjugate to $H$, then we identify $(C_{H'},p_{H'})$ with $(C_H,p_H)$ since there exists an isometry $f\: C_H\rightarrow C_{H'}$ such that $p_H=p_{H'}\circ f$.
Thus, we regard the pair $(C_H,p_H)$ as a geometric object corresponding to the conjugacy class $[H]$.
In particular, if $H$ is cyclic, then $C_H$ is homeomorphic to a circle, and $(C_H,p_H)$ can be considered as a closed geodesic on $\gS$.

Assume that $H$ is non-cyclic. Then, $C_H$ is a hyperbolic surface of finite area with geodesic boundary, whose boundary components are closed geodesics. Viewed in relation to $\gS$, the pair $(C_H,p_H)$ can be regarded as an \emph{isometrically immersed hyperbolic surface with geodesic boundary in $\gS$} since $p_H$ is the restriction of the covering map $\HH/ H\to \gS$. Note that $C_H$ has no boundary if and only if $H$ is a finite-index subgroup of $\pi_1(\gS)$.

In general, any isometrically immersed hyperbolic surface $(S,s)$ with geodesic boundary in $\gS$, where $S$ is a hyperbolic surface of finite area with closed geodesic boundary and $s\: S\to \gS$ is an isometric immersion, is isometric to a convex core $(C_H,p_H)$ for some finitely generated subgroup $H$ of $\pi_1(\gS)$; that is, there exists an isometry $f \: S \to C_H$ such that $s = p_H \circ f$. This implies that (the isometric equivalence class of) such an isometrically immersed hyperbolic surface corresponds to a conjugacy class $[H]$ via the convex core $(C_H,p_H)$.

Let $\gG=\Map (\gS)$ for simplicity.
As a consequence of the above correspondence, the quantity
\[
\# \{ [H'] \in \Gamma ([H] ) \mid \ell_\SC ([H'] ) \leq L \}
\]
in Theorem \ref{thmintro:1} can be interpreted as the number of convex cores (or isometrically immersed hyperbolic surfaces) of \emph{type} $ (C_H, p_H) $ whose $ \ell_\SC $-length is less than or equal to $L$.
Here, for a finitely generated subgroup $K$ of $\pi_1 (\gS)$, the pairs $(C_H,p_H)$ and $(C_K,p_K)$ are of the same type if there exists $[\phi]\in \Map (\gS)$ such that $\phi \circ p_H$ is homotopic to $p_K$, that is, there exist a homeomorphism $f\: C_H\rightarrow C_K$ and a continuous function $\Psi \: C_H\times [0,1]\rightarrow \gS$ such that $\Psi(\cdot,0)=\phi \circ p_H$ and $\Psi (\cdot ,1)=p_K\circ f$. Note that $\phi \circ p_H$ is homotopic to $p_K$ if and only if $[\phi ]([H])=[K]$.
This generalizes the notion that two closed geodesics are of the same type, which is used in \cite{Mir08,Mir16,ES22}.

The notion of \emph{same type} can be naturally generalized to \emph{same type with respect to $\gG$} for a proper subgroup $\gG$ of $\Map (\gS)$, by replacing $[\phi]\in\Map(\gS)$ with $[\phi]\in\gG$ in the above definition.

When $H$ is non-cyclic and $p_H$ is an embedding, $C_H$ is identified with $p_H(C_H )$ and considered as a subsurface of $\gS$. Then, the boundary $\partial C_H$ of $C_H$ is a simple multicurve on $\gS$, and we have
\[
\# \{ [H'] \in \Gamma ([H] ) \mid \ell_\SC ([H'] ) \leq L \} =s_\gG ([H]) \# \left\{ \gamma \in \Gamma (\partial C_H ) \; \middle| \; \frac{1}{2}\ell (\gamma  ) \leq L \right\},
\]
where $s_\gG ([H])$ is either $1$ or $2$ (see Lemma \ref{lem:B(H)/H is finite}). Hence, Theorem \ref{thmintro:1} follows directly from Mirzakhani's results.
The fraction $\frac{1}{2}$ in the right-hand side comes from the continuity of the `boundary projection' $\B$ (see Subsection \ref{sec:idea of Thmintro1} and \ref{sec:gen of Thmintro 1}).

We remark on the last assertion of Theorem \ref{thmintro:1}.
If $H_1,\dots,H_m$ are cyclic, then the formal weighted sum $a_1[H_1]+\cdots+a_m[H_m]$ represents a weighted multicurve, and hence Theorem \ref{thmintro:1} recovers Mirzakhani's result on counting weighted multicurves.
More generally, such a formal sum can be interpreted as a finite union of isometrically immersed hyperbolic surfaces or subsurfaces.
For example, suppose $H_1$ and $H_2$ are non-cyclic, $p_{H_1}$ and $p_{H_2}$ are embeddings, and $p_{H_1}(C_{H_1})\cap p_{H_2}(C_{H_2})=\emptyset$. In this case, $[H_1]+[H_2]$ corresponds to the union of the subsurfaces $p_{H_1}(C_{H_1})\cup p_{H_2}(C_{H_2})$.

\subsection{Boundary projection and generalized length functional}\label{sec:idea of Thmintro1}
The boundary projection $\B$, which plays an essential role in Theorem \ref{thmintro:1}, is defined as follows.
For a non-trivial, non-cyclic, finitely generated subgroup $H$ of $\pi_1(\gS)$, define
\[ \B ([H]) = \frac{1}{2}\sum_{\text{$c\:$boundary component of $C_H$}} c, \tag{$\ast$} \]
where each boundary component $c$ of $C_H$ is regarded as a closed geodesic on $\gS$ via the projection $p_H\: C_H\to \gS$.
Hence, $\B([H])$ is a weighted multicurve on $\gS$.

Note that if $H$ is a finite-index subgroup of $\pi_1(\gS)$, then $C_H$ is a hyperbolic surface without boundary; that is, $\partial C_H$ is empty and hence $\B([H]) = 0$. From the viewpoint of counting, the orbit $\Map (\gS) ([H])$ is a finite set, and hence the left-hand side of the formula in Theorem \ref{thmintro:1} is zero.

In the case where $H$ is cyclic, $C_H$ itself is a closed geodesic on $\gS$, and we define $\B ([H])=C_H$.
Specifically, in this context, we can interpret the boundary components of $C_H$ as the set consisting of two formal copies of $C_H$.
Indeed, if we associate a surface to $H$, then it is an annulus each of whose boundary components is homotopic to $C_H$.
This interpretation enables us to use the above definition $(\ast)$ for a cyclic subgroup $H$ of $\pi_1(\gS)$.

The generalized length functional $\ell_\SC$ in Theorem \ref{thmintro:1} is defined as the composition $\ell \circ \B$.
More generally, by using the result of \cite{ES22}, $\ell_\SC$ can be replaced with $F\circ \B$ for any positive homogeneous continuous function $F$ on $\GC_K(\gS)$. 
In addition, if $\B([H])\not=0$, that is, $H$ is not a finite-index subgroup of $\pi_1(\gS)$, then $\mathfrak{c}^\Gamma_{g,r}(H)$ is a certain positive integral multiple of $\mathfrak{c}^\Gamma_{g,r}(\B([H]))>0$ (see Theorem \ref{thm:counting subgroup use length}).

Theorem \ref{thmintro:1} can be directly derived from Mirzakhani's result and Lemma \ref{lem:B(H)/H is finite}, which states that the map $\B\: \gG ([H])\rightarrow \GC (\gS)$ is finite-to-one.
This approach is similar to the one used for counting arcs in \cite{Bel23} (see \cite[Theorem 1.1 and Corollary 3.6]{Bel23} for details).

\subsection{Generalization of Theorem \ref{thmintro:1}}\label{sec:gen of Thmintro 1}
After proving Theorem \ref{thmintro:1} (Theorem \ref{thm:counting subgroup use length}) in Section \ref{sec:counting subgroup}, we will extend our discussion to a more general asymptotic formula for counting subgroups in Section \ref{sec:general counting}.
Since this generalization draws inspiration from the work presented in \cite{ES22}, we will first review the underlying principles of their proof.
Note that the space $\GC(\gS)$ of geodesic currents on $\gS$ serves as a measure-theoretic completion of the set of weighted multicurves on $\gS$; in particular it includes $\ML (\gS)$.
For $L\geq 0$, any weighted multicurve $\gamma_0$ and any finite-index subgroup $\gG$ of $\Map (\gS)$, we can define the counting measure on $\GC(\gS)$ as
\[ m_{\gamma_0}^L=\frac{1}{L^{6g-6+2r}}\sum_{\gamma \in \gG (\gamma_0 )} \delta_{\frac{1}{L}\gamma}, \]
where $\delta_{\frac{1}{L}\gamma}$ represents the Dirac measure at $\frac{1}{L}\gamma \in \GC (\gS )$.

Erlandsson-Souto \cite{ES22} proved that when $L$ goes to $\infty$, the measure $m_{\gamma_0}^L$ converges to $\mathfrak{c}^\Gamma_{g,r}(\gamma_0) m_\Thu $ with respect to the weak-$\ast$ topology on the space of Radon measures on $\GC_K(\gS)$, where $K\subset \gS$ is any compact subset including all closed geodesics in $\Gamma (\gamma_0)$.
This implies that when $L$ goes to $\infty$, $m_{\gamma_0}^L(F^{-1}([0,1]))$ converges to $\mathfrak{c}^\Gamma_{g,r}(\gamma_0) m_\Thu (F^{-1}([0,1]))$. Note that
\begin{align*}
m_{\gamma_0}^L(F^{-1}([0,1]))=
&\frac{1}{L^{6g-6+2r}}\# \left\{ \gamma\in \Gamma (\gamma_0 )\; \middle| \; F\left(\frac{1}{L}\gamma \right)\leq 1\right\} \\
=&\frac{1}{L^{6g-6+2r}}\#\{ \gamma\in \Gamma (\gamma_0 )\mid F(\gamma )\leq L\} .
\end{align*}
We consider the convergence of the measure $m_{\gamma_0}^L$ to $\mathfrak{c}^\Gamma_{g,r}(\gamma_0) m_\Thu$ as the essence of the counting problem independent of the measurement $F$. Figure \ref{fig:GC_conv} at the beginning of Section \ref{sec:general counting} illustrates the convergence of the counting measure $m_{\gamma_0}^L$.

In the case of counting subgroups, the space $\SC(\gS)$ of subset currents on $\gS$, introduced in \cite{KN13}, plays the same role as $\GC(\gS)$.
We usually assume that the limit set of a finitely generated subgroup $H$ contains at least two points.
From the viewpoint of the counting, if the limit set of $H$ has only one point, then $\gS$ has cusps and a generator of $H$ is peripheral, which implies that the $\gG$ orbits of $[H]$ is at most finite.

For any finitely generated subgroup $H$ of $\pi_1(\gS)$ whose limit set has at least two points, we can define a corresponding subset current $\eta_H\in \SC (\gS)$.
We need to keep in mind that the correspondence between the conjugacy class $[H]$ and $\eta_H$ is finite-to-one (see Proposition~\ref{prop:property of eta} and the subsequent discussion).
Nonetheless, this does not introduce any significant issues for counting subgroups.
By the correspondence and the denseness of the set
\[\{ c\eta_H\mid c>0, H\: \text{finitely generated subgroup of }\pi_1(\gS) \} \]
in $\SC(\gS)$, we can consider the space $\SC (\gS)$ as a measure-theoretic completion of weighted (sum of) conjugacy classes of finitely generated subgroups of $\pi_1(\gS)$.
The action of $\Map (\gS)$ on the set of conjugacy classes of finitely generated subgroups of $\pi_1(\gS)$ can extend to the continuous action on $\SC(\gS)$.
Note that if $H$ is cyclic, then $\eta_H$ can be regarded as a geodesic current. Hence, $\SC(\gS)$ includes $\GC(\gS)$ as a closed subspace and also includes $\ML (\gS)$.

Geometrically, $\eta_H$ is associated with the convex core $C_H$, serving as an extension of the relationship between a geodesic current and a closed geodesic on $\gS$.
When $\gS$ is closed, there are many interesting continuous functionals on $\SC(\gS)$ coming from the geometric structure of $C_H$ (see \cite{Sas22} for detail).
For example, there exists a unique continuous $\RRR$-linear functional $\Area\: \SC(\gS) \rightarrow \RRR$ such that for any non-trivial finitely generated subgroup $H$ of $\pi_1(\gS)$ we have
\[ \Area (\eta_H )=\text{area of }C_H.\]
This is indicative of the space $\SC (\gS)$ being an effective completion.

In the context of this paper, the continuous extension $\B \: \SC (\gS)\rightarrow \GC (\gS)$ of the boundary projection $\B$ is fundamental. The fraction $\frac{1}{2}$ in the definition of $\B$ is critical for the continuous extension of $\B$ to have the property that the restriction of $\B$ to $\GC (\gS)$ is the identity map.

When $\gS$ has cusps, these cusps constitute obstructions to the continuous extension of certain geometric invariants of $C_H$ (see \cite[Section 6]{Sas22b}).
Nevertheless, by restricting the domain to $\GC_K (\gS)$ for any compact subset $K$ of $\gS$, we are able to construct continuous functionals with this domain.
For example, the (geometric) intersection number $i$ of closed geodesics can be extended to a continuous functional
\[ i\: \GC_K (\gS )\times \GC (\gS )\to \RRR \]
but cannot be extended to a continuous functional on $\GC (\gS )\times \GC (\gS)$ (see \cite[Section 6]{Sas22b}). In Subsection \ref{subsec:area functional} we are going to construct the continuous area functional $\Area$ on $\SC_K(\gS):=\B^{-1}(\GC_K (\gS))$, which will be used for the proof of Theorem \ref{thmintro:2} presented below.

The weighted finite sum $a_1[H_1]+\cdots +a_m[H_m]$ in Theorem \ref{thmintro:1} corresponds to 
\[ \eta=a_1\eta_{H_1}+\cdots +a_m\eta_{H_m} \in \SC (\gS ).\]
We assume that $\B (\eta)\not=0$. Then, similarly as above, we can define the counting measure $m_\eta^L$ on $\SC (\gS)$ as
\[ m_\eta^L=\frac{1}{L^{6g-6+2r}}\sum_{\eta' \in \gG (\eta )} \delta_{\frac{1}{L}\eta '}. \]
Then, our second result can be stated as follows:
\begin{theoremintro}[See Theorem \ref{thm:general counting thm}]\label{thmintro:2}
There exists a positive integer $s_\gG (\eta)$ such that when $L$ goes to $\infty$, the counting measure $m_\eta^L$ converges to $s_\gG( \eta) \mathfrak{c}^\Gamma_{g,r}(\B(\eta )) m_\Thu $ with respect to the weak-$\ast$ topology on the space of
Radon measures on $\SC_K(\gS)$, where $K\subset \gS$ is any compact subset including all closed geodesics in $\Gamma (\B (\eta))$.
\end{theoremintro}

See Figure \ref{fig:SC_conv}, which appears before Theorem \ref{thm:general counting thm}, illustrating the convergence of $m_\eta^L$ in this theorem.

As a corollary to Theorem \ref{thmintro:2}, we have
\begin{corintro}[See Corollary \ref{cor:counting eta}]\label{corintro3}
For any positive homogeneous continuous function $F\: \SC_K(\gS) \rightarrow \RRR$, we have
\[ \lim_{L\rightarrow \infty}\frac{\# \{ \eta' \in \Gamma (\eta )\mid F (\eta )\leq L\}}{L^{6g-6+2r}}=s_\gG (\eta ) \mathfrak{c}^\Gamma_{g,r}(\B (\eta )) m_\Thu (F^{-1}([0,1])).\]
Furthermore, $F$ needs only to be positive on $\GC_K(\gS)$, not necessarily over all of $\SC_K(\gS)$. 
\end{corintro}
The last assertion can be deduced from the fact that the functional $\Area$ is positive on $\SC_K(\gS)\setminus \GC_K (\gS)$ and maintains $\Map (\gS )$-invariance (see Remark \ref{rem:counting and area} for detail). The length functional $\ell_\SC =\ell \circ \B$ satisfies the condition of $F$.
With a small modification, $\eta$ can be replaced with the formal weighted sum $a_1[H_1]+\cdots +a_m[H_m]$ in Corollary \ref{corintro3} (see Corollary \ref{cor:count general subgroups}).

\subsection{Acknowledgments}
I am sincerely grateful to Prof. Katsuhiko Matsuzaki for many enlightening conversations and his invaluable feedback.  
I would also like to thank the referee for their careful reading of the manuscript and for providing valuable comments.

This work was supported by 
Grant-in-Aid for JSPS Fellows 21J01271 and
JSPS KAKENHI Grant Number JP25K06985.

\section{Preliminaries}

In this introductory section, we establish the fundamental notations and concepts that will be used throughout this paper.

Let $\HH$ represent the hyperbolic plane. We denote the group of orientation-preserving isometries of $\HH$ by $\Isom ( \HH )$. A \emph{hyperbolic surface} is a quotient space $\HH/G$ for a torsion-free discrete subgroup $G$ of $\Isom (\HH )$.
We identify the fundamental group $\pi_1(\HH/ G)$ of the hyperbolic surface $\HH/ G$ with the subgroup $G$.
The canonical projection from $\HH$ to $\HH / G$ is denoted by $\pi$.

In this paper, we focus on hyperbolic surfaces with finite area.
A hyperbolic surface $\gS$ with finite area is either an orientable closed surface of genus $g \geq 2$ or an orientable surface of genus $g$ with $r$ cusps satisfying the condition that $2-2g+r<0$.
The latter surface is referred to as a cusped hyperbolic surface with finite area.
Here, we exclude the special case where $(g,r)=(0,3)$ since the mapping class group of such a hyperbolic surface is finite.

The \emph{limit set} of a discrete subgroup $G$ of $\Isom (\HH)$, denoted by $\gL (G)$, is the set of accumulation points of the orbit $G(x)$ in the (ideal) boundary $\partial \HH$ for $x\in \HH$, which is independent of the choice of $x$.
Note that the hyperbolic surface $\HH/ G$ has finite area if and only if $G$ is finitely generated and $\gL (G)=\HH$.

\subsection{Subset current, Geodesic current, and Measured lamination}

In this subsection, we introduce subset currents on a hyperbolic surface.  
For fundamental results on subset currents on closed and cusped hyperbolic surfaces, we refer the reader to \cite{Sas22} and \cite{Sas22b}, respectively.  
These two references will be used frequently throughout the paper.
The paper \cite{KN13}, which originally introduced the notion of subset currents, provides detailed motivation and foundational properties of the theory; see also \cite{Sas15} for the study of subset currents on free groups. For geodesic currents on closed hyperbolic surfaces, see the foundational works \cite{Bon86,Bon88}.
Moreover, for geodesic currents on cusped hyperbolic surfaces, \cite{BIPP21} provides useful background and context.

Consider the boundary $\partial \HH$ of the hyperbolic plane $\HH$. We define the space
\[ \H (\partial \HH )= \{ S\subset \partial \HH \mid \# S \geq 2\text{ and }S \text{ is compact}\} ,\]
endowed with the Vietoris topology. Note that the Vietoris topology is equivalent to the topology induced by the Hausdorff distance on $\partial \HH$, with respect to a compatible metric.
The space $\H (\partial \HH )$ is a locally compact separable metrizable space (see \cite[Theorem 2.2]{Sas22}).
The continuous action of $\Isom (\HH )$ on $\partial \HH$ extends naturally to a continuous action on $\H (\partial \HH )$.

For any set $S\in \H (\partial \HH )$, its \emph{convex hull} $\CH (S)$ is defined as the smallest convex subset of $\HH$ including all bi-infinite geodesics in $\HH$ that connect pairs of points in $S$.
This concept provides a geometric interpretation of a point of $\H (\partial \HH)$.
For example, if $S=\{ \alpha, \beta\}$, then $\CH (S)$ is the bi-infinite geodesic connecting $\alpha$ to $\beta$.

\begin{definition}[Subset current]
Let $\gS$ be a hyperbolic surface and let $G$ be the fundamental group of $\gS$. Assume that $\gL (G) =\HH$.
Note that $G$ acts on $\H (\partial \HH )$.
A \emph{subset current} on $\gS$ is defined as a $G$-invariant, locally finite (positive) Borel measure on $\H (\partial \HH )$.
A Borel measure $\mu$ is said to be \emph{locally finite} if $\mu (K)<\infty$ for any compact set $K$.
Since $\H (\partial \HH )$ is a locally compact separable metrizable space, any locally finite Borel measure on $\H (\partial \HH )$ is regular (see \cite[2.18 Theorem]{Rud86}) and thus qualifies as a Radon measure.

The space of subset currents on $\gS$, denoted by $\SC (\gS)$, is equipped with an $\RRR$-linear structure.
Furthermore, we endow $\SC (\gS)$ with the weak-$\ast$ topology.
A sequence $\{ \mu _n\}$ of $\SC (\gS)$ converges to $\mu\in \SC (\gS)$ if and only if for any compactly supported continuous function $f\: \H (\partial \HH )\rightarrow \RR$ we have
\[ \lim_{n\rightarrow \infty }\int f d\mu_n = \int f d \mu .\]
We note that $\SC (\gS)$ is a locally compact, separable, and completely metrizable space (see \cite[Theorem 2.6]{Sas22} for the case of a closed hyperbolic surface and \cite[Proposition 3.4.]{Sas22b}) for the case of a cusped hyperbolic surface).
\end{definition}

\begin{definition}[Geodesic current and Measured lamination]
A \emph{geodesic current} on a hyperbolic surface $\gS$ is a subset current whose support is included in
\[ \partial_2 \HH =\{ S \subset \partial \HH \mid \# S =2 \} . \]
Recall that the support of a measure $\mu$, denoted by $\supp (\mu )$, is the smallest closed subset whose complement has measure zero with respect to $\mu$.
The space of geodesic currents on $\gS$, denoted by $\GC (\gS)$, inherits the subspace topology from $\SC (\gS)$.

A \emph{measured lamination} is a geodesic current satisfying the condition that for any $S_1,S_2$ in its support with $S_1\not=S_2$, their convex hulls $\CH (S_1)$ and $\CH (S_2)$ are disjoint.
The space of measured laminations on $\gS$, denoted by $\ML (\gS)$, also inherits the subspace topology from $\SC (\gS)$.
Note that $\GC (\gS)$ is a closed subspace of $\SC (\gS)$, and $\ML (\gS)$ is a closed subspace of $\GC (\gS)$.

For any compact subset $K$ of $\gS$, a geodesic current $\mu$ is said to be supported on $K$ if for any $S\in \supp (\mu)$ the convex hull $\CH (S)$ is included in $\pi^{-1} (K)\subset \HH$.
Let $\GC_K(\gS)$ denote the set of geodesic currents supported on $K$. This is a closed subspace of $\GC (\gS)$.
In the case where $\gS $ is a cusped hyperbolic surface, we often focus on the subspace $\GC_K(\gS)$ rather than  on $\GC (\gS)$. This approach, which is standard in the study of measured laminations and geodesic currents on cusped hyperbolic surfaces, has also been employed in the previously cited references \cite{ES22, Sas22b, BIPP21}.
When $\gS$ is a closed hyperbolic surface, we usually take $K = \gS$, so that $\GC_K(\gS) = \GC(\gS)$.
\end{definition}

Let $\gS$ be a hyperbolic surface of finite area and let $G$ be its fundamental group.
We denote by $\Sub (G)$ the set of finitely generated subgroups of $G$ whose limit sets have at least two points, which ensures that $\gL (H)\in \H (\partial \HH )$ for any $H\in \Sub (G)$.
We remark that the limit set $\gL(H)$ has only one point if and only if $\gS$ is a cusped hyperbolic surface and $H$ is generated by a single element whose representative is peripheral.

For $H\in \Sub (G)$ we define a Borel measure $\eta_H$ on $\H (\partial \HH)$ by
\[ \eta_H=\sum_{gH \in G/H } \delta _{g\gL (H)}, \]
where $\delta_{g\gL(H)}$ is the Dirac measure at $g\gL (H)$.
It is straightforward to verify that $\eta_H$ is $G$-invariant.
Furthermore, we can show that $\eta_H$ is locally finite (see \cite[Lemma 2.7]{Sas22} and \cite[Theorem 2.6 and Remark 2.7]{Sas22b}).
Hence, $\eta_H$ qualifies as a subset current on $\gS$.

If $H$ is a cyclic subgroup generated by $\gamma \in G$, then $\eta_H$ is a geodesic current and the above construction fits with the construction of a geodesic current associated to a closed geodesic freely homotopic to a representative of $\gamma$. Moreover, if $\gamma$ has a simple representative, then $\eta_H$ is a measured lamination.

The current $\eta_H$ associated to $H$ has the following properties:

\begin{proposition}[See {\cite[Proposition 2.10]{Sas22}}]\label{prop:property of eta}
For $H_1, H_2\in \Sub (G)$,
\begin{enumerate}
\item if $H_1$ is a $k$-index subgroup of $H_2$, then $\eta_{H_1}=k\eta_{H_2}$;
\item if $H_1$ is conjugate to $H_2$, then $\eta_{H_1}=\eta_{H_2}$.
\end{enumerate}

\end{proposition}

Based on the property (2) from the above proposition, we consider that $\eta_H$ corresponds to a conjugacy class $[H]$ for $H\in \Sub (G)$.
We remark that this correspondence is not one-to-one but one-to-finite.

To elaborate, there may exist $H_1,H_2\in \Sub (G)$ such that $H_1$ is not conjugate to $H_2$ but $\eta_{H_1}=\eta_{H_2}$.
In such a case, $\gL(H_1)=\gL(H_2)$ up to conjugacy. We then observe that both $H_1$ and $H_2$ are $k$-index subgroups of 
the stabilizer
\[ \Stab (\gL (H_1) ):= \{ g\in G \mid g (\gL (H_1) )=\gL (H_1) \}\]
for some $k\geq 2$.
The point is that the number of $k$-index subgroups of $\Stab (\gL (H_1))$ is finite.
As a result, the number of conjugacy classes $[H]$ satisfying the condition $\eta_H=\eta_{H_1}$ is finite.

In general, for any positive integer $k$, the number of $k$-index subgroups of a finitely generated group $J$ is finite. This is because a $k$-index subgroup of $J$ appears as the stabilizer for some action of $J$ on a set of $k$-elements, which is determined by the action of the finite generators of $J$ on this finite set.

For $H\in \Sub (G)$, we consider the convex hull $\CH (\gL(H))$ of the limit set $\gL(H)$, which is the smallest convex $H$-invariant subset of $\HH$ including all bi-infinite geodesics connecting pairs of points in $\gL (H)$.
Then, the \emph{convex core} $C_H$ of $H$ is defined as the quotient space $\mathrm{CH}(\gL (H))/H$. The universal covering map $\pi\: \HH \rightarrow \gS$ induces the projection $p_H\: C_H\rightarrow \gS$. Note that the convex core $C_H$ can be identified with the smallest closed convex subset of the quotient space $\HH/H$ such that the inclusion map is a homotopy equivalence. Then, the map $p_H$ is the restriction of the covering map $\HH/ H\rightarrow \gS$ to $C_H$.

We remark that the convex core $C_H$ can be considered as a geometric object corresponding to $\eta_H$.
In particular, if $H$ is a cyclic subgroup generated by $\gamma \in G$, then the projection $p_H \: C_H\rightarrow \gS$ represents an unoriented closed geodesic freely homotopic to a representative of $\gamma$.
In the context of geodesic currents, the geodesic current $\eta_{\langle \gamma \rangle}$ is usually identified with the closed geodesic $c$ freely homotopic to a representative of $\gamma$.
When $H$ is not cyclic, $C_H$ becomes a hyperbolic surface of finite area with closed geodesic boundary, which is referred to as a \emph{hyperbolic surface of finite type}. If the projection $p_H$ is injective, then we can identify $C_H$ with the subsurface $p_H(C_H)$ of $\gS$.

Note that we do not assume closed geodesics to be primitive. For $\gamma \in G$ whose representative is freely homotopic to a primitive closed geodesic, we can consider the geodesic current $\eta_{\langle \gamma^k\rangle }$ for any positive integer $k$. Then, we have
\[ \eta_{\langle \gamma^k\rangle }=k\eta_{\langle \gamma \rangle }\]
since $\langle \gamma^k\rangle$ is a $k$-index subgroup of $\langle \gamma \rangle$ by the property (1) of Proposition \ref{prop:property of eta}.

\subsection{Action of mapping class group on subset currents}

Let's consider an orientation-preserving homeomorphism $\phi \: \gS \rightarrow \gS$. We can take a lift $\tilde{\phi}\: \HH \rightarrow \HH$ of $\phi$ such that $\pi\circ \tilde{\phi}=\phi \circ \pi$. The lift $\tilde{\phi}$ is a homeomorphism satisfying the condition $\tilde{\phi}G \tilde{\phi}^{-1}=G$.
Furthermore, $\tilde{\phi}$ induces the self-homeomorphism $\partial \tilde{\phi}\: \partial \HH \rightarrow \partial \HH$, and moreover, induces the self-homeomorphism $\partial \tilde{\phi}\: \H (\partial \HH )\rightarrow \H (\partial \HH)$.

When we take a subset current $\mu\in \SC (\gS)$, then we can define $\phi(\mu ) \in \SC(\gS)$ as the push-forward measure $(\partial \tilde{\phi})_\ast(\mu)$ by $\partial \tilde{\phi}$, explicitly,
\[ \phi (\mu )(E)=\mu ((\partial \tilde{\phi})^{-1}(E))\]
for any Borel subset $E\in \H (\partial \HH)$.

For an orientation-preserving homeomorphism $\psi$ isotopic to $\phi$ and its lift $\tilde{\psi}\: \HH \rightarrow \HH$, there exists $g\in G$ such that $g\circ \partial \tilde{\phi}=\partial \tilde{\psi}$. As $\mu\in \SC (\gS)$ is $G$-invariant, we get $\phi(\mu)=\psi(\mu)$. Consequently, the action of the mapping class group 
\[ \Map (\gS):= \{ \phi \mid \phi \: \gS \rightarrow \gS,\text{ an orientation-preserving homeomorphism} \} /\mathrm{isotopy}\]
on $\SC(\gS)$ is established. Both $\ML(\gS)$ and $\GC(\gS)$ are $\Map(\gS)$-invariant subsets of $\SC(\gS)$.

For $[\phi_0] \in \Map(\gS)$, we can choose a representative $\phi$ of $[\phi_0]$ such that $\phi$ fixes the base point of the fundamental group $G$. Then, we can consider $\phi$ as an automorphism of $G$, and we have
\[ \partial \tilde{\phi}(\gL (H))= \gL (\phi (H))\]
for an appropriate lift $\tilde{\phi}$ of $\phi$ and any $H\in \Sub (G)$.
Hence, 
$[\phi_0](\eta_H)=\eta_{\phi (H)}$.

Note that the subset current $\eta_H$ is determined by the conjugacy class $[H]$ of $H$, and $[\phi ]\in \Map(\gS)$ can be considered as an outer automorphism of $G$.
Unless there are any issue, we simplify $[\phi] \in \Map (\gS)$ to $\phi$.
With this convention, for $\phi \in \Map (\gS)$, we can express
\[ \phi (\eta_H)=\eta_{\phi(H)},\]
which allows us to consider the action of $\Map (\gS)$ on $\SC(\gS)$ as a continuous extension of the action of $\Map(\gS)$ on the set of conjugacy classes $\Sub(G)/{\sim}$.
Furthermore, the action of $\Map (\gS)$ on $\GC (\gS)$ is a continuous extension of its action on the set of all closed geodesics on $\gS$.

Note that if $H_1,H_2\in \Sub(G)$ are $k$-index subgroups of $H\in \Sub(G)$, then
\[ \eta_{H_1}=k\eta_H =\eta_{H_2}.\]
In such a case, $H_1$ is not necessarily conjugate to $H_2$, and there may exist $\phi\in \Map(\gS)$ such that $\phi([H_1])=[H_2]$ (see the following example).

\begin{example}\label{example:free group}
Consider $x,y\in G=\pi_1(\gS)$ as shown in Figure \ref{fig:genus_2_surface}.
The subgroup $\langle x,y \rangle$ generated by $x,y$ is isomorphic to the free group of rank $2$. We identify the convex core $C_{\langle x,y \rangle}$ with the subsurface $p_{\langle x,y \rangle}(C_{\langle x,y\rangle})$ of $\gS$ since $p_{\langle x,y \rangle}\: C_{\langle x,y\rangle}\rightarrow \gS$ is injective.

We regard $\langle x,y\rangle$ as the fundamental group of the wedge of two circles, denoted by $R_2$.
Put
\[ H=\langle x^4, xy,y^2, x^2yx, x^2y^{-1}x\rangle ,\]
which corresponds to the covering graph $\Delta_H$ of $R_2$ as shown in the left side of Figure \ref{fig:Delta_H_graph}.
The covering graph allows us to conclude that $H$ is a $4$-index subgroup of $\langle x,y\rangle$.

Let $\phi$ be the Dehn twist about $x$ fixing the base point of $G$ such that $\phi (x)=x, \phi (y)=xy$.
Then,
\[ \phi (H)=\langle x^4, x^2y, xy xy, x^3yx, x^2y^{-1}\rangle,\]
which corresponds to the covering graph $\Delta_{\phi(H)}$ of $R_2$ as shown in the right of Figure \ref{fig:Delta_H_graph}.
It is clear that $H$ is not conjugate to $\phi(H)$ in $\langle x,y\rangle$ as $\Delta_H$ is not isomorphic to $\Delta_{\phi(H)}$.

We note that the injectivity of $p_{\langle x,y\rangle}$ implies that for any $g\in G\setminus \langle x,y \rangle$ we have
\[ g\CH (\gL (\langle x,y\rangle ))\cap \CH ((\langle x,y \rangle))=\emptyset ,\]
and so $g\langle x,y\rangle g^{-1}\cap \langle x,y \rangle= \{ \id \}$. As a result, we conclude that $H$ is not conjugate to $\phi(H)$ in $G$.
Therefore, for the conjugacy class $[H]$, we have $\phi([H])\not=[H]$ and
\[ \eta_{\phi(H)}=4\eta_{\langle x,y \rangle }=\eta_{H}.\]

\begin{figure}[h]
\centering
\includegraphics{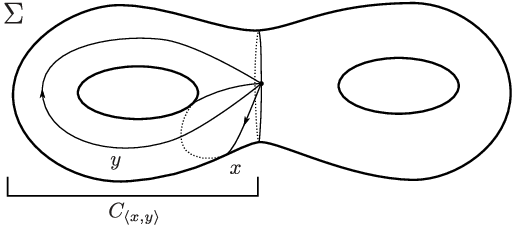}
\caption{The convex core $C_{\langle x,y\rangle}$ is described as the subsurface of $\gS$.}\label{fig:genus_2_surface}
\end{figure}
\begin{figure}[h]
\centering
\includegraphics{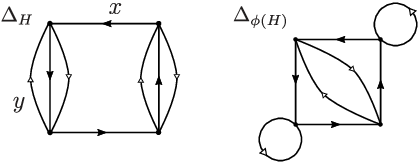}
\caption{The left of the figure is the covering graph $\Delta_H$ and the right of the figure is the covering graph $\Delta_{\phi(H)}$.}\label{fig:Delta_H_graph}
\end{figure}
\end{example}

From the above example, we see that the stabilizer
\[ \mathrm{Stab}(\eta_H):=\{\phi \in \Map (\gS) \mid \phi (\eta_H)=\eta_H\}\]
of $\eta_H$ includes $\mathrm{Stab}([H])$ as a possibly proper subset. Furthermore, the proposition below is presented. Notably, this proposition can also be independently derived as a corollary of Lemma \ref{lem:B(H)/H is finite}.

\begin{proposition}\label{lem:stab(eta_H)/stab([H])}
For any $H\in \Sub (G)$, the index $[ \Stab(\eta_H )\: \Stab([H])]$ is finite.
\end{proposition}
\begin{proof}
Let $H_0=\Stab_G (\gL (H))=\{ g\in G\mid g\gL (H)=\gL(H)\}$, which includes $H$.
Then, we have $\gL(H)=\gL (H_0)$ and both $C_H$ and $C_{H_0}$ are hyperbolic surfaces of finite type.
Hence, the canonical projection $p_H^{H_0}\:C_H \rightarrow C_{H_0}$ is a finite-covering, which implies that $H_0$ includes $H$ as a finite-index subgroup. Let $k$ be the index of $H$ in $H_0$.
Therefore, $k\eta_{H_0}=\eta_H$.

For $K\in \Sub (G)$, if $\eta_K=\eta_{H_0}$, then there exists $g\in G$ such that $g\gL(K)=\gL(H_0)$.
By the definition of $H_0$, $gKg^{-1}$ is included in $H_0$. Since 
\[ \eta_{gKg^{-1}}=\eta_K =\eta_{H_0},\]
$gKg^{-1}$ is a $1$-index subgroup of $H_0$, which implies that $gKg^{-1}=H_0$.
As a result, $\eta_K=\eta_{H_0}$ if and only if $K$ is conjugate to $H_0$.
Hence, $\Stab (\eta_{H_0} )=\Stab ([H_0])$.

Take any $[\phi] \in \Stab(\eta_H)$. Then, $[\phi]$ fixes $\eta_{H_0}=\frac{1}{k}\eta_H$, and consequently $[\phi]$ fixes $[H_0]$.
By considering $[\phi]$ as an outer automorphism of $G$, we see that $[\phi]([H])$ equals $[J]$ for a $k$-index subgroup $J$ of $gH_0g^{-1}$ for some $g\in G$.
Hence, we obtain the action of $\Stab (\eta_H)$ on the set
\begin{align*}
A
&=\{ [J]  \mid J\text{ is a $k$-index subgroup of $gH_0 g^{-1}$ for some $g\in G$} \}\\
&=\{ [J] \mid J\text{ is a $k$-index subgroup of $H_0$}\}.
\end{align*}

Recall that the number of $k$-index subgroups of a given finitely generated subgroup is finite, which implies that $A$ is a finite set.
The action of $\Stab (\eta_H )$ on $A$ induces the one-to-one correspondence between $\Stab(\eta_H)/\Stab ([H])$ and the orbit $\Stab(\eta_H) ([H]) \subset A$.
As a result, the index $[ \Stab(\eta_H )\: \Stab([H])]$ is bounded above by the cardinality $\# A$.
\end{proof}

\subsection{Counting curves and Thurston measure}

Let $\gS$ be a hyperbolic surface of genus $g$ with $r$ cusps.
We allow $r$ to be $0$, and then $\gS$ is a closed hyperbolic surface of genus $g\geq 2$.
We assume that $\gS$ is not a thrice-punctured sphere since the mapping class group of a thrice-punctured sphere is a finite group.

A closed curve is said to be \emph{essential} if it is neither null-homotopic nor peripheral. Note that we do not assume closed curves to be primitive.
Since our surface $\gS$ has a hyperbolic structure, a free homotopy class of an essential closed curve contains a unique closed geodesic.
Recall that for a cyclic subgroup $\langle \gamma \rangle \in \Sub (G)$, a geodesic current $\eta_{\langle \gamma \rangle }$ corresponds to an unoriented closed geodesic $c$ on $\gS$ that is free homotopic to a representative of $\gamma$. In this case, $\eta_{\langle \gamma \rangle }$ can be referred as $\eta_c$.
Closed geodesics on $\gS$ that we consider are usually unoriented and regarded as geodesic currents.

We define a \emph{weighted multicurve} $\gamma$ of $\gS$ as a weighted finite sum
\[ a_1\gamma_1+\cdots + a_m \gamma_m ,\]
where $a_1,\dots ,a_m\in \RRR$ and $\gamma_1,\dots ,\gamma_m$ are closed geodesics.
A weighted multicurve $\gamma$ is regarded as a geodesic current, that is,
\[ \gamma = a_1\eta_{\gamma_1}+\cdots +a_m\eta_{\gamma_m}.\]
If $\gamma_1,\dots ,\gamma_m$ are pair-wise disjoint simple closed geodesics, then $\gamma$ is called a \emph{weighted simple multicurve} and considered as an element of $\ML(\gS)$. Note that a simple closed geodesic is primitive.

For a closed geodesic $\gamma$, we denote by $\ell(\gamma)$ the length of $\gamma$.
From \cite[Theorem 3.9]{ES22}, $\ell$ can be continuously extended to an $\RRR$-linear functional
\[ \ell \: \GC _K(\gS ) \rightarrow \RRR\]
for any compact subset $K$ of $\gS$, that is, for any closed geodesic $\gamma$, we have
\[ \ell (\eta_\gamma )=\ell ( \gamma ).\]
Moreover, for a weighted multicurve $a_1\gamma_1+\cdots + a_m \gamma_m\in \GC _K(\gS )$, we have
\[ \ell (a_1\gamma_1+\cdots + a_m \gamma_m )=a_1\ell (\gamma_1 )+\cdots +a_m \ell (\gamma _m).\]

\begin{definition}[Thurston measure]
Let $\ML_\ZZ(\gS)$ be the subset of $\ML (\gS)$ consisting of simple weighted multicurves whose weights are positive integer.
The \emph{Thurston measure} $m_\Thu$ on $\ML(\gS)$ is defined as the weak-$\ast$ limit of 
\[ \frac{1}{L^{6g-6+2r}}\sum_{\gamma \in \ML_\ZZ(\gS)}\delta_{\frac{1}{L}\gamma}\]
on $\ML(\gS)$ when $L$ goes to $\infty$ (see \cite[Theorem 4.16]{ES22}).
We consider $m_\Thu$ as the measure on $\SC(\gS)$ (or on $\GC(\gS)$) by defining
\[ m_\Thu (\SC(\gS) \setminus \ML (\gS) )=0.\]
\end{definition}

Note that $6g-6+2r$ is the dimension of $\ML (\gS)$, and $m_\Thu$ belongs to the Lebesgue measure class on $\ML(\gS)$.
The convergence in the above definition is an analogue of the following convergence in $\RR^n$:
\[ \lim_{L\rightarrow \infty} \frac{1}{L^n}\sum_{x\in \ZZ^n}\delta_{\frac{1}{L}x} =\text{Lebesgue measure on } \RR^n .\]

Now, we state Mirzakhani's curve counting theorem:

\begin{theorem}[\cite{Mir08,Mir16}]\label{thm:Mirzakhani}
For any weighted multicurve $\gamma_0$ and any finite-index subgroup $\Gamma$ of $\Map(\gS)$, there exists a positive constant $\mathfrak{c}^\Gamma_{g,r}(\gamma_0)$ such that
\[ \lim_{L\rightarrow \infty}\frac{\# \{ \gamma \in \Gamma (\gamma_0 )\mid \ell (\gamma )\leq L\}}{L^{6g-6+2r}}=\mathfrak{c}^\Gamma_{g,r}(\gamma_0) m_\Thu (\ell^{-1}([0,1])).\]
\end{theorem}

Erlandsson-Souto \cite{ES22} extended this theorem into a general version where the length functional $\ell$ can be replaced with any positive homogeneous continuous function $F$ on $\GC_K(\gS)$ for a compact subset $K$ of $\gS$ that includes all closed geodesics in $\Gamma (\gamma_0)$. A function $F\: \GC_K(\gS)\rightarrow \RRR$ is \emph{positive} if $F(\mu)>0$ for any non-zero $\mu\in \GC_K(\gS)$, and $F$ is \emph{homogeneous} if $F(c\mu)=cF(\mu)$ for any $c\in \RRR$ and $\mu\in \GC_K(\gS)$.
This result is an immediate corollary of the following theorem, which can be regarded as fundamental to the counting problem.

\begin{theorem}[{\cite[Theorem 8.1]{ES22}}]\label{thm:ES}
For any weighted multicurve $\gamma_0$ and any finite-index subgroup $\Gamma$ of $\Map(\gS)$, there exists a positive constant $\mathfrak{c}^\Gamma_{g,r}(\gamma)$ such that
\[ \lim_{L\rightarrow \infty}\frac{1}{L^{6g-6+2r}}\sum_{\gamma \in \Gamma(\gamma_0)}\delta_{\frac{1}{L}\gamma}=\mathfrak{c}^\Gamma_{g,r}(\gamma_0 ) m_\Thu. \]
Here the convergence takes place with respect to the weak-$\ast$ topology on the space of Radon measures on $\GC_K(\gS)$, where $K\subset \gS$ is any compact subset including all closed geodesics in $\Gamma (\gamma_0)$.
\end{theorem}

For the convenience of later reference, we states the corollary explicitly:

\begin{corollary}[{\cite[Theorem 9.1]{ES22}}]\label{cor:ES}
In the setting of the above theorem, for any positive homogeneous continuous function $F\: \GC_K(\gS)\rightarrow \RRR$, we have
\[ \lim_{L\rightarrow \infty}\frac{\# \{ \gamma \in \Gamma (\gamma_0 )\mid F (\gamma )\leq L\}}{L^{6g-6+2r}}=\mathfrak{c}^\Gamma_{g,r}(\gamma_0) m_\Thu (F^{-1}([0,1])).\]
\end{corollary}

The outline of the proof of the corollary is as follows.
From the assumption, we can obtain the following convergence of the measures of the set $F^{-1}([0,1])$:
\[ \lim_{L\rightarrow \infty}\left(\frac{1}{L^{6g-6+2r}}\sum_{\gamma \in \Gamma(\gamma_0)}\delta_{\frac{1}{L}\gamma}\right) (F^{-1}([0,1]))=\mathfrak{c}^\Gamma_{g,r}(\gamma_0 ) m_\Thu (F^{-1}([0,1])).\]
The measure in the limit of the left-hand side is equal to
\begin{align*}
&\frac{1}{L^{6g-6+2r}}\# \left\{ \gamma\in \Gamma (\gamma_0 )\; \middle| \; F\left(\frac{1}{L}\gamma \right)\leq 1\right\} \\
=&\frac{1}{L^{6g-6+2r}}\#\{ \gamma\in \Gamma (\gamma_0 )\mid F(\gamma )\leq L\} .
\end{align*}
Hence, we obtain the general version of Theorem \ref{thm:Mirzakhani}.
Note that the case of weighted \emph{simple} multicurves of Theorem \ref{thm:ES} is first proved by Mirzakhani in \cite{Mir08}.

\section{Counting subgroups}\label{sec:counting subgroup}

Let $\gS$ be a hyperbolic surface of genus $g$ with $r$ cusps and let $G$ be the fundamental group of $\gS$.
We allow $r$ to be $0$, and then $\gS$ is a closed hyperbolic surface of genus $g\geq 2$.
We assume that $\gS$ is not a thrice-punctured sphere.
Recall that $\Sub (G)/{\sim}$ is the set of conjugacy classes of finitely generated subgroups of $G$ whose limit sets have at least two points.

Our goal in this section is to extend the ``asymptotic formula'' for weighted multicurves in Theorem \ref{thm:Mirzakhani} to the asymptotic formula for $\eta_H\in \SC (\gS)$ associated with $[H]\in \Sub(G)/{\sim}$.
As a corollary, we also derive the asymptotic formula for $[H]\in \Sub(G)/{\sim}$ itself.

When considering the counting problem for the finitely generated subgroups of $G$, one challenge is determining how to measure their `length'.
The continuous projection $\B\: \SC (\gS)\rightarrow \GC (\gS)$, introduced in \cite[Chapter 7]{Sas22}, plays a fundamental role in this measurement. In fact, utilizing the projection $\B$ along with Lemma \ref{lem:B(H)/H is finite in general}, we can derive the asymptotic formula for the conjugacy class of a finitely generated subgroup.

We will review the construction of $\B$ in Subsection \ref{subsec:conti proj}. Here, we overview the property of $\B$ from the geometric viewpoint.
Recall that if $H\in \Sub (G)$ is non-cyclic, then the convex core $C_H$ becomes a hyperbolic surface of finite area possibly with geodesic boundary, referred to as a hyperbolic surface of finite type.
We use $\partial C_H$ to denote the set of connected components of the boundary of $C_H$.

In the case where $H$ is cyclic, $C_H$ is homeomorphic to a circle. Specifically, in this context, we define $\partial C_H$ as the set consisting of two formal copies of $C_H$, expressed as:
\[ \partial C_H=\{ C_H, C_H'\}.\]
The point is that $\eta_{C_H}=\eta_{C_{H'}}=\eta_H$. Remark that if we associate a surface to $H$, then it is an annulus whose boundaries are both homotopic to $C_H$.

With the above notation, for any $H\in \Sub (G)$ we have
\[ \B (\eta_H )=\frac{1}{2}\sum_{c\in \partial C_H}\eta_c,\]
where $\eta_c$ is a geodesic current on $\gS$ corresponding to the closed geodesic $p_H(c)$.
In essence, the map $\B$ is the $\RRR$-linear continuous extension of the mapping from $C_H$ to half the sum of its boundary components.
Such an extension is unique since the set 
\[ \{ c\eta_H \mid c>0,\ H\in \Sub (G) \}\]
of \emph{rational subset currents} is dense in $\SC(\gS)$ (see \cite[Theorem 8.21]{Sas22} and \cite[Theorem 2.13]{Sas22b}).

Note that if $H$ is cyclic, then $\mathcal{B}(\eta_H)=\eta_H$ by the definition of $\partial C_H$.
If $H$ is a finite-index subgroup of $G$, then $\partial C_H$ is empty and $\B (\eta_H)=0$, the zero measure.

\subsection{Continuous projection from subset currents to geodesic currents}\label{subsec:conti proj}

In \cite{Sas22}, the author introduced a continuous projection, denoted by $\mathcal{B}$, that maps from the space $M(\H (\partial \HH) )$ of locally finite Borel measures on $\H (\partial \HH)$ to its subspace $M(\partial_2\HH)$.
Herein, we review the construction of $\mathcal{B}$.

Let $\mathcal{O}$ be the set of all open intervals of $\partial \HH$.
Let $\mathcal{M}$ be the counting measure on $\mathcal{O}$, that is, 
\[ \mathcal{M}(U)=(\text{the cardinality of }U) \in \RRR \cup \{ \infty \}.\]
for any $U\subset \mathcal{O}$.
For any $S\in \H(\partial \HH)$ we define $b(S)$ as the subset of $\mathcal{O}$ consisting of all connected components of $\partial \HH \setminus S$. 
We then define a map $\varphi \: \H(\partial \HH)\times \mathcal{O}\rightarrow \{ 0,1\}$ as
\[ \varphi( S, \alpha )=\delta_\alpha (b(S)),\]
that is, $\varphi(S,\alpha)=1$ if $\alpha \in b(S)$; $\varphi(S,\alpha)=0$ if $\alpha \not\in b(S)$.

For any $\mu\in M (\H (\partial \HH ))$ the measure $\mathcal{B}(\mu)$ on $\partial_2\HH$ is defined as follows: for any Borel subset $E\subset \partial_2\HH$
\[ \mathcal{B}(\mu) (E)=\frac{1}{2}\int_{\H (\partial \HH )}\left( \int_{b(E)} \varphi (S,\alpha )d\mathcal{M} (\alpha )\right) d\mu (S),\]
where $b(E):=\cup_{S\in E}b(S)$, a subset of $\mathcal{O}$. Note that for any Borel subsets $E_1,E_2\subset \partial_2\HH$, $E_1\cap E_2=\emptyset$ if and only if $b(E_1)\cap b(E_2)=\emptyset$.

The strategy for the proof of the continuity of $\B$ with respect to the weak-$\ast$ topology on $M (\H (\partial \HH ))$ is outlined as follows (detailed in the proof of \cite[Theorem 7.1]{Sas22}). Consider any $\mu \in M (\H (\partial \HH ))$ and a sequence $\{ \mu_n \} \subset M (\H (\partial \HH ))$ converging to $\mu$. 
According to \cite[Proposition 5.45]{Sas22}, known as the Portmanteau theorem, it is sufficient to prove that for any relatively compact subset $E\subset \partial_2\HH$ with $\mu (\partial E)=0$, the following limit holds:
\[ \lim_{n\rightarrow \infty} \B (\mu_n )(E)=\B (\mu ) (E).\]
This is established by showing that the map
\[ \H (\partial \HH ) \ni S\mapsto \int_{b(E)}  \varphi (S,\alpha )d\mathcal{M} (\alpha )\]
is a bounded function on $\H (\partial \HH )$ with compact support and the set of its non-continuous points has measure zero with respect to $\mu$. Utilizing \cite[Proposition 5.45]{Sas22} once more confirms that $\lim_{n\rightarrow \infty} \B (\mu_n )(E)=\B (\mu ) (E)$.

If the support of $\mu\in M(\H (\partial \HH ))$ is included in $\partial_2\HH$, then for any Borel subset $E\subset \partial_2\HH$ and any $S\in \partial_2\HH$ we have
\[ \int_{b(E)}  \varphi (S,\alpha )d\mathcal{M} (\alpha )= \begin{cases}
2 & (S \in E )\\
0 & (S \not\in E)
\end{cases}.
\]
Hence,
\begin{align*}
\mathcal{B}(\mu )(E)&=\frac{1}{2}\int_{\H (\partial \HH )}\left( \int_{b(E)}  \varphi (S,\alpha )d\mathcal{M} (\alpha )\right) d\mu (S)\\
&=\frac{1}{2}\int_{\partial_2 \HH}\left( \int_{b(E)}  \varphi (S,\alpha )d\mathcal{M} (\alpha )\right) d\mu (S)\\
&=\frac{1}{2}\int_E 2 d\mu (S)\\
&=\mu(E).
\end{align*}
This implies that the restriction of $\mathcal{B}$ to $M(\partial_2\HH)$ is the identity mapping.
We remark that the fraction $\frac{1}{2}$ is crucial to derive this property.

Now, we consider a self-homeomorphism $f\: \partial \HH \rightarrow \partial \HH$, which induces the self-homeomorphism $f$ of $\partial_2\HH$ and of $\H (\partial \HH)$.
For any $\mu \in M (\H (\partial \HH ))$ and any Borel subset $E\subset \partial_2\HH$ we have
\begin{align*}
\mathcal{B}(f_\ast (\mu ))(E)&=\frac{1}{2}\int_{\H (\partial \HH )}\left( \int_{b(E)}  \varphi (S,\alpha )d\mathcal{M} (\alpha )\right) df_\ast(\mu )(S)\\
&=\frac{1}{2}\int_{\H (\partial \HH )}\left( \int_{b(E)}  \varphi (f(S),\alpha )d\mathcal{M} (\alpha )\right) d\mu (S)\\
&=\frac{1}{2}\int_{\H (\partial \HH )}\left( \int_{b(E)}  \varphi (S,f^{-1}(\alpha ) )d\mathcal{M} (\alpha )\right) d\mu (S)\\
&=\frac{1}{2}\int_{\H (\partial \HH )}\left( \int_{b(f^{-1}(E))}  \varphi (S,\alpha )d\mathcal{M} (\alpha )\right) d\mu (S)\\
&=\mathcal{B}(\mu )(f^{-1}(E))\\
&=f_\ast (\mathcal{B}(\mu ))(E).
\end{align*}
From the above we derive the following lemma.

\begin{lemma}
Let $J$ be a group that acts on $\partial \HH$ continuously. If $\mu$ is a $J$-invariant locally finite Borel measure on $\H (\partial \HH)$, then so is $\mathcal{B}(\mu)$.
In particular, for a hyperbolic surface $\gS$ of finite area, $\B$ is a continuous $\RRR$-linear map from $\SC(\gS)$ to $\GC(\gS)$ whose restriction to $\GC (\gS)$ is the identity map. Additionally, $\mathcal{B}\: \SC(\gS)\rightarrow \GC (\gS)$ is a $\Map (\gS)$-equivariant map, that is,
\[ \mathcal{B}(\phi (\mu ))=\phi (\mathcal{B}(\mu ))\]
for any $\phi \in \Map (\gS)$ and any $\mu\in \SC (\gS)$.
\end{lemma}

Finally, we check that for any $H\in \Sub (G)$ the following equality holds:
\[ \B (\eta_H )=\frac{1}{2}\sum_{c\in \partial C_H}\eta_c .\]
For any Borel subset $E\subset \partial_2\HH$,
\begin{align*}
&2\B (\eta_H )(E)\\
=&\int_{\H (\partial \HH )} \left( \int_{b(E)} \varphi (S ,\alpha ) d\mathcal{M}(\alpha ) \right) d\eta_H (S)\\
=&\sum_{gH\in G/H} \int_{b(E)} \varphi (g\gL(H) ,\alpha ) d\mathcal{M}(\alpha )
=\sum_{gH\in G/H} \int_{b(E)} \delta_{\alpha }(b(g\gL( H))) d\mathcal{M}(\alpha ) \\
=&\sum_{gH\in G/H} \int_{b(g \gL (H))} \delta_{\alpha }(b(E)) d\mathcal{M}(\alpha )
=\sum_{gH\in G/H} \sum_{\alpha \in b(\gL (H))}\delta_{g(\alpha ) }(b(E)) \\
=&\sum_{gH\in G/H} \sum_{\alpha \in b(\gL (H) )}g_\ast (\delta_{ \partial \alpha})(E)
=\sum_{gH\in G/H} g_\ast \left( \sum_{\alpha \in b(\gL (H) )}\delta_{ \partial \alpha}\right) (E)\\
\stackrel{\star}{=}&\sum_{gH\in G/H} g_\ast \left( \sum_{c\in \partial C_H} \sum_{h\langle c\rangle \in H/\langle c\rangle } \delta_{ h\gL (\langle c \rangle )}  \right) (E)\\
=&\sum_{c\in \partial C_H} \sum_{gH\in G/H} g_\ast \left( \sum_{h\langle c\rangle \in H/\langle c\rangle } \delta_{ h\gL (\langle c \rangle )}  \right) (E)\\
=&\sum_{c\in \partial C_H} \sum_{g\langle c \rangle \in G/\langle c \rangle} \delta_{ g\gL (\langle c \rangle )}  (E)
=\sum_{c\in \partial C_H} \eta_c (E).
\end{align*}
See \cite[Lemma 7.2]{Sas22} for the detail of the equality $\stackrel{\star}{=}$.
In the above calculation, the closed geodesic $c\in \partial C_H$ also represents an element of $G=\pi_1(\gS)$ whose representative is free homotopic to $p_H(C_H)$.

\subsection{Counting subgroups via boundary lengths of convex cores}

Recall that the length function $\ell$ of $\gS$ can be continuously extended to an $\RRR$-linear functional
\[ \ell \: \GC _K(\gS ) \rightarrow \RRR,\]
for any compact subset $K$ of $\gS$.
For the preimage $\SC_K(\gS)=\B^{-1}(\GC_K (\gS))$, we define $\ell_\SC$ as
\[ \ell_\SC=\ell\circ \B \: \SC_K(\gS )\rightarrow \RRR .\]
We can express
\[ \ell_\SC (\eta_H )=\frac{1}{2}\sum_{c\in \partial C_H}\ell (c) \]
for $H\in \Sub (G)$. For any $H_1,\dots ,H_m\in \Sub (G)$ and $a_1,\dots ,a_m\in \RRR$ we define
\[ \ell_\SC (a_1[H_1]+\cdots +a_m[H_m]):=\ell_\SC (a_1\eta_{H_1}+\cdots +a_m\eta_{H_m})=\sum_{i=1}^m a_i \ell_\SC (\eta_{H_i} ).\]
We will use the `length' $\ell_\SC$ for counting `subgroups'.

Recall that $\B\: \SC (\gS )\rightarrow \GC (\gS)$ is a $\Map(\gS)$-equivariant map. For any $H\in \Sub (G)$ and any finite-index subgroup $\gG$ of $\Map (\gS)$ the surjective map
\[ \B|_{\gG (\eta_H) }\: \gG (\eta_H )\rightarrow \gG (\B (\eta_H ))\]
is not injective in general. However, we observe that $\B|_{\gG (\eta_H )}$ is a finite-to-1 map when $\B (\eta_H)\not=0$.
Recall that when $H$ is a finite-index subgroup of $G$, then $C_H$ has no boundary.
In this case, the orbit $\Map (\gS) ([H])$ is included in the finite set
\[ \{ [J] \in \Sub(G)/{\sim} \mid [G: J]=[G: H]\} .\]

The following lemma plays a fundamental role in proving the asymptotic formula for $\eta_H$.

\begin{lemma}\label{lem:B(H)/H is finite}
For any finite-index subgroup $\gG$ of $\Map (\gS)$ and any $H\in \Sub (G)$ with $\B(\eta_H)\not=0$, we define
\[ s_\gG ([ H] ):=\# \Stab_\gG (\B(\eta_H ))/ \Stab_\gG ([H] )\]
where $\# \Stab_\gG (\B(\eta_H ))=\gG \cap \Stab (\B (\eta_H))$ and $\Stab_\gG ([H] )=\gG \cap \Stab([H] )$.
Then, $s_\gG([H])$ is finite.
In addition, 
\[ s_\gG(\eta_H ):= \# \Stab_\gG (\B(\eta_H ))/ \Stab_\gG (\eta_H )\]
is also finite.
Consequently, $\B_{\gG (\eta_H )}$ is $s_\gG (\eta_H )$-to-$1$, that is, for any $\gamma \in \gG (\B (\eta_H )) $ we have
\[ s_\gG (\eta_H )=\# (\B |_{\gG (\eta_H )} )^{-1}(\gamma ).\]
\end{lemma}
\begin{proof}
We remark that $\Stab_\gG (\eta_H)$ includes $\Stab_\gG ([H])$. Hence, it is sufficient to prove that $s_\gG ([H])$ is finite.

First, consider the case where the canonical projection $p_H\: C_H\rightarrow \gS$, induced by the canonical projection $\pi \:\HH \rightarrow \gS$, is injective.
In this case, we can regard $C_H$ as a subsurface of $\gS$, where each boundary component is a simple closed geodesic on $\gS$.
Conversely, any such subsurface of $\gS$ induces a conjugacy class of a finitely generated subgroup of $G$, which coincides with $[H]$ in this instance.
Note that if we have another injective projection $p_J\: C_J \rightarrow \gS$ for $J\in \Sub (G)$ such that $p_J(C_J)=p_H(C_H)$, then $p_J^{-1}\circ p_H\: C_H\rightarrow C_J$ induces a covering isomorphism
\[ \HH /H\rightarrow \HH /J,\]
which implies that $H$ is conjugate to $J$.

Observe that $p_H$ is injective if and only if for any $g\in G\setminus H$ we have
\[ g\CH ( \gL (H))\cap \CH ( \gL (H))=\emptyset.\]
When this condition is satisfied, for any orientation-preserving homeomorphism $\phi\: \gS\rightarrow \gS$ fixing the base point of $G=\pi_1(\gS)$ and its lift $\tilde{\phi}\: \HH\rightarrow \HH$ we have
\[ \tilde{\phi}(g\CH ( \gL (H)))\cap \tilde{\phi}(\CH ( \gL (H)))=\emptyset.\]
Hence,
\[ \phi(g) \CH (\gL (\phi(H) ))\cap \CH (\gL (\phi(H) ))=\emptyset.\]
Therefore, $p_{\phi(H)}\: C_{\phi(H)}\rightarrow \gS$ is also injective.
This means that the action of $\Map (\gS)$ on $\Sub (G)/{\sim}$ preserves the injectivity of the canonical projection $p_H$.

From the above discussions, we see that $\Stab_\gG (\B (\eta_H ))$ acts on the set
\[ \{ [J] \in \Sub (G)/{\sim} \mid  p_J\: C_J\rightarrow \gS \text{ is injective} \text{ and }\B (\eta_J )=\B(\eta_H )\} ,\]
which consists of two elements: one is $[H]$ and the other corresponds to the complementary subsurface $\gS\setminus C_H$.
As a result, $s_\gG ([H])$ is either $1$ or $2$.
In fact, $s_\gG(\eta_H )=s_\gG([H])=2$ if $C_H$ is homeomorphic to the closure of $\gS\setminus C_H$ (see Figure \ref{fig:genus_2_surface} in Example \ref{example:free group}).

Next, let's consider the general case where $p_H\: C_H\rightarrow \gS$ is not necessarily injective. By the main result in \cite{Sco78,Sco85}, there exists a finite-index subgroup $G_0$ of $G$ such that $H$ is a subgroup of $G_0$ and the canonical projection $p_H^{G_0}\: C_H\rightarrow C_{G_0}$ is injective.
Then, $C_H$ can be regarded as a subsurface of $C_{G_0}$. Each component $c$ of $\partial C_H$ is a simple closed geodesic of $C_{G_0}$. Since $p_{G_0}\:C_{G_0}\rightarrow \gS$ is a finite-covering, the lifts of $p_H(c)=p_{G_0}(p_H^{G_0}(c))$ to $C_{G_0}$ consist of a finite collection of geodesics.
Let $k$ be the index of $G_0$ in $G$.

Similarly to the above case, we have the action of $\Stab_\gG (\B (\eta_H) )$ on the set
\[  A=\left\{ [J]\in \Sub (G)/{\sim} \; \middle|
\begin{array}{l}
G'\in \Sub (G), J\subset G'\subset G, [G:G']=k, \\
p_J^{G'}\:C_J\rightarrow C_{G'}\text{ is injective, and } \B (\eta _J) =\B (\eta _H)
\end{array}
\right\} .\]
Since $\B (\eta _J) =\B (\eta _H)$, each boundary component $c_0$ of $C_J$ must coincide with one of the lifts of $p_H(c)$ to $C_{G'}$ for some boundary component $c$ of $C_H$.
This implies that the number of boundary component candidates for the subsurface $C_J$ of $C_{G'}$ is finite.
Noting that the number of $k$-index subgroups of a finitely generated group is finite, we can conclude that $A$ is a finite set. Therefore, the cardinality $s_\gG ([H])$, which equals the cardinality of the orbit $\Stab_\gG (\B (\eta_H ))([H])$ within $A$, is finite.
\end{proof}

Under the conditions of the above lemma, we have
\[ \Stab_\gG(\B (\eta_H ))/\Stab_\gG ([H] )\cong \Stab_\gG(\B (\eta_H ))/\Stab_\gG (\eta_H ) \times \Stab_\gG (\eta_H )/\Stab_\gG ([H]),\]
which implies that the index $[\Stab_\gG (\eta_H ): \Stab_\gG ([H])]$ is also finite. Hence, we obtain Proposition \ref{lem:stab(eta_H)/stab([H])} as a corollary.

We can extend the above lemma to the situation of weighted sum of conjugacy classes of finitely generated subgroups of $G$.

\begin{lemma}\label{lem:B(H)/H is finite in general}
Let $\gG$ be a finite-index subgroup of $\Map (\gS)$. Given any $H_1,\dots ,H_m\in \Sub (G)$, consider the formal weighted sum
\[ J=a_1[H_1]+\cdots +a_m[H_m] \quad (a_1,\dots ,a_m\in \RRR)\]
and the weighted sum
\[ \eta =a_1\eta_{H_1}+\cdots +a_m\eta_{H_m} .\]
Assume that $\B (\eta)\not=0$. Then, both cardinalities
\[ s_\gG ( J ):= \# \Stab_\gG (\B (\eta ))/ \Stab_\gG (J ) \text{ and } s_\gG ( \eta ):=\# \Stab_\gG (\B(\eta ))/ \Stab_\gG (\eta )\]
are finite. Consequently, $\# \Stab_\gG (\eta )/\Stab_\gG (J )$ is also finite.
\end{lemma}

\begin{proof}
Note that the action of $\Map (\gS)$ on $\Sub (G)/{\sim}$ extends linearly to the set of weighted sums of $\Sub(G)/{\sim}$.

Our approach to this lemma follows the same pattern as the preceding lemma.
First, we observe that $\Stab_\gG (\eta)$ includes $\Stab_\gG ( J)$, and that $\Stab_\gG (J)$ includes
\[ \Stab_\gG ([H_1],\dots ,[H_m]):=\{ \phi \in \gG \mid \phi ([H_i]) =[H_i]\ (i=1,\dots ,m)\} \]
as a finite-index subgroup. This is established by considering the canonical group homomorphism from $\Stab_\gG (J)$ to the permutation group of $m$ elements.
To complete the proof, it is sufficient to prove that
\[ \# \Stab_\gG (\B(\eta ))/\Stab_\gG ([H_1],\dots ,[H_m])\]
is finite.

Recall that for $H\in \Sub (G)$, $\partial C_H$ is the set of boundary components of $C_H$, and can be regarded as a subset of $\GC (\gS)$. We introduce $\partial \eta$ as:
\[ \partial \eta=\bigcup_{i=1}^m \partial C_{H_i} \subset \GC (\gS ) . \]
Then, $\Stab_\gG (\B(\eta ))$ acts on $\partial \eta$, which is a finite set.
The point is that for each $i=1,\dots ,m$ and for any $\phi \in \Stab_\gG (\B (\eta ))$, $\phi ( \B (\eta_{H_i}) )$ is not necessarily equal to $\B (\eta_{H_i})$ but $\phi (\partial C_{H_i})$ is included in $\partial \eta$.

Similarly to the proof of Lemma \ref{lem:B(H)/H is finite}, for each $i=1,\dots ,m$, we can take a finite-index subgroup $G_i$ of $G$ such that $G_i$ includes $H_i$ and
$p_{H_i}^{G_i} \: C_{H_i}\rightarrow C_{G_i}$
is injective. Let $k_i$ denote the index of $G_i$ in $G$.
Then, we see that $\Stab_\gG (\B(\eta ))$ acts on the finite set
\[  A_i=\left\{ [H]\in \Sub (G)/{\sim} \; \middle|
 \begin{array}{l}
G'\in \Sub (G), H\subset G'\subset G, [G:G']=k_i, \\
p_H^{G'}\:C_H\rightarrow C_{G'}\text{ is injective, and } \partial C_H \subset \partial \eta
 \end{array}
\right\}. \]
Moreover, we have the diagonal action of $\Stab_\gG (\B (\eta ))$ on the finite set
\[ A_1\times \cdots \times A_m,\]
which contains $([H_1], \dots ,[H_m])$. Therefore, $\# \Stab_\gG (\B(\eta ))/\Stab_\gG ([H_1],\dots ,[H_m])$ is finite.
\end{proof}

As a corollary to the above lemma, we immediately obtain the following theorem.

\begin{theorem}\label{thm:counting subgroup use length}
Let $\gS$ be a hyperbolic surface of genus $g$ with $r$ cusps and assume that $(g,r)\not=(0,3)$.
Let $\gG$ be a finite-index subgroup of $\Map (\gS)$. Given any $H_1,\dots ,H_m\in \Sub (G)$, consider
\[ J=a_1[H_1]+\cdots +a_m[H_m] \text{ and } \eta =a_1\eta_{H_1}+\cdots +a_m\eta_{H_m}\quad  (a_1,\dots ,a_m\in \RRR) .\]
Assume that $\B (\eta)\not=0$. Then, we have 
\[ \lim_{L\rightarrow \infty}\frac{\# \{ J' \in \Gamma (J )\mid \ell_\SC (J' )\leq L\}}{L^{6g-6+2r}}=s_\gG (J )\mathfrak{c}^\Gamma_{g,r}(\B (\eta )) m_\Thu (\ell^{-1}([0,1]))\]
and
\[ \lim_{L\rightarrow \infty}\frac{\# \{ \eta' \in \Gamma (\eta )\mid \ell_\SC (\eta' )\leq L\}}{L^{6g-6+2r}}=s_\gG (\eta )\mathfrak{c}^\Gamma_{g,r}(\B (\eta )) m_\Thu (\ell^{-1}([0,1])),\]
where the constant $\mathfrak{c}^\Gamma_{g,r}(\B (\eta ))$ originates from Theorem \ref{thm:Mirzakhani}, and the constants $s_\gG (J )$, $s_\gG (\eta )$ come from Lemma \ref{lem:B(H)/H is finite in general}.
Moreover, $\ell_\SC$ can be replaced with $F\circ \B$ for any positive homogeneous continuous function $F\: \GC_K (\gS) \rightarrow \RRR$.
\end{theorem}
\begin{proof}
From Lemma \ref{lem:B(H)/H is finite in general}, we have
\begin{align*}
&\#\{ J' \in \Gamma (J )\mid \ell_\SC (J' )\leq L\} \\
=&\#\{  \phi \in \Gamma/\Stab_\gG (J) \mid \ell_\SC (\phi (J) )\leq L\} \\
=&\#\{  \phi \in \Gamma/\Stab_\gG (J) \mid \ell (\B (\phi (\eta )) )\leq L\} \\
=&\#\{  \phi \in \Gamma/\Stab_\gG (J) \mid \ell (\phi (\B (\eta )) )\leq L\} \\
=&\#\{  (\phi, \psi ) \in \Gamma/\Stab_\gG (\B (\eta ))\times \Stab_\gG (\B (\eta ))/\Stab_\gG (J) \mid \ell (\phi (\psi (\B (\eta ))) )\leq L\} \\
=&s_\gG (J) \#\{ \phi \in \Gamma/\Stab_\gG (\B (\eta ))\mid \ell (\phi (\B (\eta )) )\leq L\} \\
=&s_\gG (J) \#\{ \gamma \in \Gamma (\B (\eta ))\mid \ell (\gamma )\leq L\}.
\end{align*}
Hence, by Theorem \ref{thm:Mirzakhani}, we have
\begin{align*}
&\lim_{L\rightarrow \infty}\frac{\# \{ J' \in \Gamma (J )\mid \ell_\SC (J' )\leq L\}}{L^{6g-6+2r}}\\
=&s_\gG (J) \lim_{L\rightarrow \infty}\frac{\#\{ \gamma \in \Gamma (\B (\eta ))\mid \ell (\gamma )\leq L\}}{L^{6g-6+2r}}=s_\gG (J )\mathfrak{c}^\Gamma_{g,r}(\B (\eta )) m_\Thu (\ell^{-1}([0,1])).
\end{align*}
The second asymptotic formula in the theorem also follows by the same argument.
Moreover, using Corollary \ref{cor:ES}, even when we replace $\ell_\SC$ with $F\circ \B$ for any positive homogeneous continuous function $F\: \GC_K (\gS) \rightarrow \RRR$, we can obtain the same asymptotic formula. Note that $m_\Thu (\ell^{-1}([0,1]))$ becomes $m_\Thu (F^{-1}([0,1]))$ in this case.
\end{proof}

While it is possible to introduce a compact notation $c_{g,r}^\Gamma (\eta)$ to denote the product $s_\gG (\eta )\mathfrak{c}^\Gamma_{g,r}(\B (\eta ))$, we choose not to use this simplification. The components $s_\gG (\eta )$ and $\mathfrak{c}^\Gamma_{g,r}(\B (\eta ))$ each play a significant and distinct role in the subsequent analysis. Therefore, to maintain the clarity of their individual contributions, we will explicitly refer to each component throughout.

\section{General counting theorem of subgroups}\label{sec:general counting}

In this section, our goal is to extend Theorem \ref{thm:ES} to the case of subgroups.

\begin{figure}[h]
\centering
\includegraphics{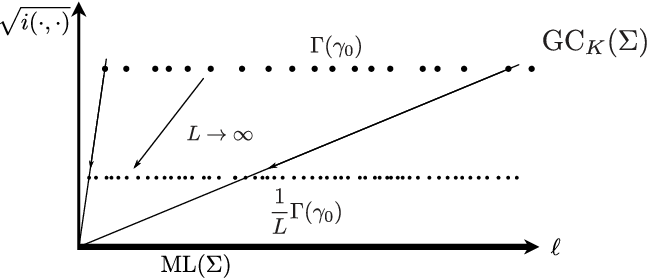}
\caption{This figure illustrates the convergence of the limit in Theorem \ref{thm:ES}.}\label{fig:GC_conv}
\end{figure}

First, we overview a certain important method used in the proof of Theorem \ref{thm:ES}.
Let $\gS$ be a hyperbolic surface of genus $g$ possibly with $r$ cusps and let $\gamma_0$ be a weighted multicurve on $\gS$.
Take a compact subset $K$ of $\gS$ including all closed geodesics in $\Map (\gS ) (\gamma_0 )$.
For any $\phi \in \Map (\gS )$, we have
\[ i\left( \frac{1}{L}\phi (\gamma_0) , \frac{1}{L}\phi (\gamma_0 )\right) =\left( \frac{1}{L}\right) ^2 i(\gamma_0 ,\gamma_0 ) \rightarrow 0 \quad (L\rightarrow \infty ).\]
Remark that the intersection number $i\: \GC_K(\gS)\times \GC_K(\gS)\rightarrow \RRR$ is a continuous $\RRR$-bilinear $\Map$-invariant functional and
\[ \ML (\gS )=\{ \mu \in \GC (\gS )\mid i(\mu ,\mu )=0 \}.\]
Then, we see that the limit of the counting measure
\[ \lim_{L\rightarrow \infty }\frac{1}{L^{6g-6+2r}}\sum_{\gamma \in \Gamma(\gamma_0)}\delta_{\frac{1}{L}\gamma}\]
in the asymptotic formula in Theorem \ref{thm:ES} is a measure supported by $\ML (\gS )$ (see Figure \ref{fig:GC_conv}).
See \cite[Proposition 6.2]{ES22} for detail.

In the following subsection, we will introduce the area functional $\Area \: \SC_K(\gS)\rightarrow\RRR$, which will play a role analogous to that of the self-intersection number $i(\cdot, \cdot )$ in proving the subset current version of Theorem \ref{thm:ES}.
Recall that $\SC_K(\gS) =\B ^{-1} (\GC_K (\gS) )$.

\subsection{Area functional}\label{subsec:area functional}

When $\gS$ is a closed hyperbolic surface, we have already obtained the area functional
\[ \Area \: \SC (\gS) \rightarrow \RRR\]
in \cite[Theorem 3.3]{Sas22}, which includes the case that $\gS$ is a higher-dimensional closed hyperbolic manifold.
The area functional $\Area$ is a continuous $\RRR$-linear functional with the property that for any $H\in \Sub (G)$, we have
\[ \Area (\eta_H )= \text{area of }C_H.\]
Moreover, we note that for any $\mu \in \SC (\gS)$, $\Area (\mu )=0$ if and only if $\mu \in \GC (\gS)$.

In this subsection, we are going to prove the existence of the area functional for hyperbolic surfaces that have cusps. Explicitly, our goal is to prove the following theorem:

\begin{theorem}\label{thm:area functional}
Let $\gS$ be a hyperbolic surface of genus $g$ with at least one cusp. There exists an $\RRR$-linear functional
\[ \Area \: \SC (\gS ) \rightarrow \RRR\]
such that for any compact subset $K$ of $\gS$, $\Area$ is continuous on $\SC_K(\gS)$, and for any $H\in \Sub (G)$ we have
\[ \Area (\eta_H )=\text{area of }C_H.\]
\end{theorem}

Let $\gS$ be a hyperbolic surface of finite area with at least one cusp. For the action of $G=\pi_1(\gS)$ on $\HH$ we take a Dirichlet fundamental domain $\F$ centered at some point. As $\gS$ is not compact, $\F$ is a non-compact finite polygon of $\HH$.
We introduce the function
\[ f_\F \: \H (\partial \HH )\rightarrow \RRR;\ S\mapsto m_\HH (CH (S) \cap \F ), \]
where $m_\HH$ represents the hyperbolic area measure on $\HH$.

Following the method used in the proof of \cite[Theorem 3.3]{Sas22}, for any $H\in \Sub (G)$ and any complete system $R$ of representatives of $G/H$, we have
\begin{align*}
\int f_\F d\eta_H 
&=\sum_{gH\in G/H} m_\HH (CH (g \gL ( H ))\cap \F )\\
&=\sum_{g\in R} m_\HH (CH ( \gL ( H ))\cap g^{-1}\F )\\
&=m_\HH \left(  \CH (\gL (H )) \cap \bigcup_{g\in R} g^{-1}\F \right)\\
&=\text{area of }C_H.
\end{align*}
The last equality holds because $T:=\CH (\gL (H )) \cap (\cup_{g\in R} g^{-1}\F)$ is a `measure-theoretic' fundamental domain for the action of $H$ on $\CH (\gL (H))$.
Explicitly, $T$ satisfies the conditions that
\[ H(T)=\CH (\gL (H ))\cap \bigcup_{g\in R} H (g^{-1} \F) =\CH (\gL (H ))\cap \HH =\CH (\gL (H)) ,\]
and that for any non-trivial $h\in H$,
\[ h(T)\cap T\subset h\left( \bigcup_{g\in R} g^{-1}\F \right) \cap \bigcup_{g\in R} g^{-1}\F = \bigcup_{g_1,g_2\in R} hg_1^{-1}\F \cap g_2^{-1}\F , \]
which is a set of measure zero with respect to $m_\HH$.

Now, we define the area functional as
\[ \Area \: \SC (\gS ) \rightarrow \RRR; \ \mu \mapsto \int f_\F d\mu.\]
It follows from the above argument that $\Area$ is an $\RR$-linear functional satisfying the condition that for any $H\in \Sub (G)$ we have
\[ \Area (\eta_H )=\text{area of }C_H.\]
However, $\Area$ is not continuous on $\SC (\gS)$ when $\gS$ has cusps.
To elaborate, according to \cite[Lemma 5.2]{Sas22b}, for two parabolic elements $\alpha, \beta \in G$ with $\alpha^\infty\not=\beta^\infty$, the sequence $\eta_{\langle \alpha ^n ,\beta^n \rangle}$ of subset currents converges to the geodesic current
\[ \eta_{ \{ \alpha ^\infty ,\beta^\infty \} }:=\sum_{g\in G } \delta_{g\{ \alpha ^\infty ,\beta ^\infty \} }, \]
where $\alpha^\infty ,\beta^\infty \in \partial \HH$ denote the fixed points of $\alpha$ and $\beta$, respectively.
Meanwhile, the area of $C_{\langle \alpha ^n ,\beta^n \rangle}$ constantly equals $2\pi$ but $\Area (\eta_{ \{ \alpha ^\infty ,\beta^\infty \} } )=0$.

We note that the function $f_\F$ is continuous on $\H (\partial \HH)$ by the proof of \cite[Proposition 3.2]{Sas22}, but the support of $f_\F$ is not compact since $\F$ is not bounded.
Hence, we can conclude that the discontinuity of $\Area$ arises from the fact that the support of $f_\F$ is not compact.

Our goal in the remaining part of this subsection is to prove that the restriction of $\Area$ to $\SC_K(\gS)$ is continuous for any compact subset $K$ of $\gS$.
In preparation for our proof, we fix a compact subset $K$ of $\gS$ and denote by $\F_K$ the intersection of the fundamental domain $\F$ and the preimage $\pi^{-1}(K)$. Note that $\F_K$ is compact.

For any subset $D\subset \HH$ we define
\[ A(D)=\{ S\in \H (\partial \HH ) \mid \CH (S) \cap D \not= \emptyset \}.\]
By \cite[Lemma 3.7, 3.8]{Sas22}, if $D$ is compact, then $A(D)$ is a compact subset of $\H (\partial \HH)$.
Hence, $A (\F_K)$ is compact.
Note that the support of $f_\F$ coincides with $A(\F)$, which is not compact.

Using the above notation, we observe that for any $\mu \in \GC_K (\gS)$, we have
\[\mu (A (\F ))=\mu (A(\F_K )).\]
To elaborate, consider any $S\in A(\F)\setminus A(\F_K)$, which means that
\[ \CH (S)\cap \F \not= \emptyset \text{ and } \CH (S)\cap (\F \cap \pi^{-1}(K))=\emptyset. \]
Then, $\CH (S)\cap (\F \setminus \pi^{-1}(K))\not= \emptyset$, indicating that $\CH(S)$ is not included in $\pi^{-1} (K)$.
By the definition of $\GC_K (\gS)$, $S$ does not belong to $\supp (\mu )$.
Therefore, 
\[ (A(\F)\setminus A(\F_K )) \cap \supp (\mu )=\emptyset,\]
which implies that $\mu (A (\F) )=\mu (A(\F_K))$.

The above property of $\GC_K (\gS)$ can be extended to the case of subset currents.

\begin{lemma}
For any $\mu \in \SC_K (\gS)$, we have
\[ \mu (A (\F ))=\mu (A (\F_K )).\]
\end{lemma}
\begin{proof}
Take any $S_0\in A (\F )\setminus A (\F_K )$. It is sufficient to see that there exists an open neighborhood $U$ of $S_0$ such that $\mu (U\cap A(\F ))=0$.
Since $S_0\in A (\F )\setminus A (\F_K )$, there exists a boundary component $c_0=\CH (\{ x_0,y_0 \} )$ of $\CH (S_0)$ such that
\[ c_0\cap \F\not= \emptyset \text{ and }c_0 \cap \F_K=\emptyset .\]
Note that if $S_0\in \partial_2\HH$, then $\{ x_0 ,y_0 \}=S_0$.
We take a small open neighborhood $U$ of $S_0$ such that $U\cap A (\F_K )=\emptyset$ and for each $S\in U$ there exists a boundary component $c=\CH (\{ x ,y\} )$ of $\CH (S)$ close to $c_0$.
Then, the set $V$ consisting of all such $\{ x, y\}$ is an open neighborhood of $\{ x_0,y_0\}$ in $\partial_2 \HH$, and $V\cap A (\F_K )=\emptyset$.

By the definition of the projection $\B \: \SC (\gS) \rightarrow \GC (\gS)$, we have
\begin{align*}
2 \B (\mu ) (V\cap A(\F ) )
=&\int_{\H (\partial \HH )} \left(\int_{b(V\cap A(\F ) )} \varphi (S, \alpha ) d\mathcal{M} (\alpha )\right) d\mu (S)\\
\geq & \int_{U\cap A(\F )  }\left( \int_{b(V\cap A(\F ) )} \varphi (S, \alpha ) d\mathcal{M} (\alpha )\right) d\mu (S).
\end{align*}
Recall that $\varphi (S, \alpha)=\delta_\alpha (b(S))=1$ if $\alpha $ is a connected component of $\partial \HH \setminus S$.
Since for any $S\in U\cap A(\F ) $ there exists $\{ x, y\} \in V\cap A(\F ) $ such that $\CH (\{ x ,y\} )$ is a boundary component of $\CH (S)$, we have $\varphi (S ,\alpha )=1$ for either interval $\alpha \in b(\{ x ,y\} )\subset b(V\cap A(\F ) )$. Hence,
\begin{align*}
2 \B (\mu ) (V \cap A(\F ) )
\geq & \int_{U \cap A(\F )  }\left(\int_{b(V\cap A(\F ) )} \varphi (S, \alpha ) d\mathcal{M} (\alpha )\right) d\mu (S)\\
\geq & \int_{U\cap A(\F )  }1 d\mu (S)=\mu (U\cap A(\F ) ).
\end{align*}
As $\B (\mu )\in \GC_K (\gS)$, we have $\B (\mu ) (V\cap A(\F ) )=0$. Hence, $\mu (U\cap A(\F ) )=0$.
\end{proof}

Now, we are going to prove Theorem \ref{thm:area functional}.

\begin{proof}[Proof of Theorem \ref{thm:area functional}]
Consider an arbitrary $\mu \in \SC_K (\gS)$ and a sequence $\{ \mu_n\}$ of $\SC_K (\gS)$ converging to $\mu$ as $n\rightarrow \infty$.
Take a compact subset $K'$ of $\gS$ such that the interior of $K'$ includes $K$, and the complement of $K'$ in $\gS$ is a union of cusp neighborhoods.
By the definition of the area functional and the previous lemma, we have
\begin{align*}
\Area (\mu )
=\int f_\F d\mu 
=&\int_{A(\F )} f_\F d\mu \\
=&\int_{A(\F_{K'} )} f_\F d\mu \\
=&\int f_\F \cdot \chi_{A (\F_{K'} )} d\mu ,
\end{align*}
where $\chi_{A(\F_{K'})}$ is the characteristic function of $A (\F_{K'} )$.
Similarly, we have
\[ \Area (\mu_n )=\int f_\F \cdot \chi_{A (\F_{K'} )} d\mu_n .\]

It is important to note that the function $f_\F \cdot \chi_{A (\F_{K'} )} $ is not continuous but its support, included in $A(\F_{K'})$, is compact.
We observe the set $\Delta (f_\F \cdot \chi_{A (\F_{K'} )} )$ of non-continuous points of $f_\F \cdot \chi_{A (\F_{K'} )} $.
Since $f_\F$ is continuous and $\F$ is a non-compact convex polygon in $\HH$, we see that for any $S\in \Delta (f_\F \cdot \chi_{A (\F_{K'} )} )$ we have
\[  \CH (S)\cap (\F \cap \pi^{-1} (\partial K' ))\not=\emptyset \text{ and } \CH(S)\cap (\F \cap \pi^{-1}(\mathrm{Int}(K') ))=\emptyset  . \]
This implies that $S$ belongs to $A(\F) \setminus A (\F_K )$ since $K$ is included in the interior $\mathrm{Int}(K')$ of $K'$.
By the previous lemma, we see that $\mu (A (\F) \setminus A (\F_K ) )=0$, which implies that $\mu (\Delta (f_\F \cdot \chi_{A (\F_{K'} )} ) )=0$.

By the Portmanteau theorem (see \cite[Proposition 5.45]{Sas22}), given that $f_\F \cdot \chi_{A (\F_{K'} )}$ is a bounded function with compact support which is $\mu$-a.e. continuous, we can conclude that $\Area(\mu_n)$ converges to $\Area (\mu)$ when $n$ tends to infinity. This completes the proof.
\end{proof}

\subsection{Counting subgroups and convergence to Thurston measure}

Let $\gG$ be a finite-index subgroup of $\Map (\gS)$.
For any multicurve $\gamma_0$, we define the counting measure
\[ m_{\gamma_0}^L= \frac{1}{L^{6g-6+2r}}\sum_{\gamma \in \Gamma (\gamma_0)} \delta_{\frac{1}{L}\gamma } \] 
on $\GC (\gS)$. From Theorem \ref{thm:ES}, it follows that when $L$ tends to $\infty$, $m_{\gamma_0}^L$ converges to a specific positive constant multiple of the Thurston measure $m_{\Thu}$, with respect to the weak-$\ast$ topology on the space of Radon measures on $\GC_K(\gS)$ for any compact subset $K$ including all closed geodesics in $\gamma_0$.
Recall that both $\GC_K(\gS)$ and $\SC_K(\gS)=\B^{-1}(\GC _K (\gS))$ are locally compact separable metrizable spaces.

Given any $H_1,\dots ,H_m\in \Sub (G)$, we consider
\[ \eta =a_1\eta_{H_1}+\cdots +a_m\eta_{H_m}\quad  (a_1,\dots ,a_m\in \RRR), \]
with the assumption $\B (\eta )\not=0$.
Let $K$ be a compact subset that includes all closed geodesics in $\gG (\B (\eta))$.
Then, we have
\[ \lim_{L\rightarrow \infty}m_{\B (\eta )}^L= \mathfrak{c}_{g,r}^\gG (\B (\eta ))m_{\Thu}.\]

Similarly to the case of multicurves, we can define the counting measure
\[ m_{\eta}^L =\frac{1}{L^{6g-6+2r}}\sum_{\eta' \in \gG (\eta )} \delta_{\frac{1}{L}\eta'} \]
on $\SC_K (\gS)$.
It is reasonable to hypothesize that when $L$ tends to $\infty$, the measure $m_{\eta}^L$ will converge to a specific positive constant multiple of $m_{\Thu}$.
Actually, we can establish the following theorem (see Figure \ref{fig:SC_conv}):

\begin{figure}[h]
\centering
\includegraphics{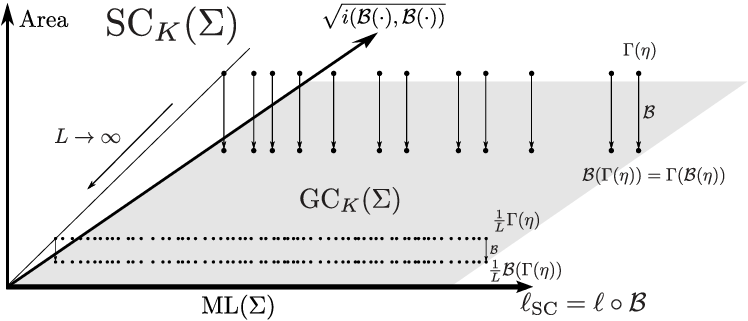}
\caption{This figure illustrates the convergence of the limit in Theorem \ref{thm:general counting thm}.
The key observation is that the area functional $\Area \: \SC_K(\gS) \to \RRR$ is a continuous $\RRR$-linear $\Map$-invariant functional, and the function $\sqrt{i(\B (\cdot ),\B (\cdot ))}\: \SC_K (\gS) \times \SC_K (\gS) \to \RRR$ is a continuous $\RRR$-linear $\Map$-invariant functional.}\label{fig:SC_conv}
\end{figure}

\begin{theorem}\label{thm:general counting thm}
Under the conditions described above, we have
\[ \lim_{L\rightarrow \infty}m_{\eta }^L= s_\gG(\eta )\mathfrak{c}_{g,r}^\gG (\B (\eta ))m_{\Thu},\]
where the convergence takes place with respect to the weak-$\ast$ topology on the space of Radon measures on $\SC_K(\gS)$.
The constant $\mathfrak{c}^\Gamma_{g,r}(\B (\eta ))$ comes from Theorem \ref{thm:Mirzakhani}, and the constant $s_\gG (\eta )$ comes from Lemma \ref{lem:B(H)/H is finite in general}.
\end{theorem}
\begin{proof}
From Theorem \ref{thm:ES}, for any continuous function $f\: \GC_K(\gS)\rightarrow \RR$ with compact support, we have
\[ \lim_{L\rightarrow \infty }\left| \int f  dm_{\B(\eta)}^L - \int f d(\mathfrak{c}_{g,r}^\gG (\B (\eta ))m_{\Thu} )\right| =0.\]
This implies
\[ \lim_{L\rightarrow \infty }\left| \int f  d(s_\gG(\eta )m_{\B(\eta)}^L) - \int f d\big( s_\gG(\eta )\mathfrak{c}_{g,r}^\gG (\B (\eta ))m_{\Thu} \big)\right| =0.\]
We can replace $f$ with any continuous function $f\: \SC_K (\gS)\rightarrow \RR$ with compact support, as $\GC_K(\gS)$ is a closed subset of $\SC_K(\gS)$.

Thus, what we need to show is that for any such continuous function $f\: \SC_K(\gS)\rightarrow \RR$ with compact support
\[ \lim_{L\rightarrow \infty }\left| \int f  dm_{\eta}^L - \int f d(s_\gG (\eta )m_{\B (\eta )}^L )\right| =0.\]

We recall the relation
\[ \gG/\Stab_\gG (\eta )\cong \gG/ \Stab_\gG (\B (\eta ))\times \Stab_\gG (\B (\eta ))/\Stab_\gG (\eta ).\]
Since $s_\gG (\eta)$ is the cardinality of $\Stab_\gG (\B (\eta ))/\Stab_\gG (\eta )$, we have
\begin{align}
&\left| \int f  dm_{\eta}^L - \int f d(s_\gG (\eta )m_{\B (\eta )}^L )\right| \notag \\ 
=&\frac{1}{L^{6g-6+2r}} \left| \sum_{\eta '\in \gG (\eta)}f\left( \frac{1}{L}\eta '\right) - s_\gG (\eta )\sum_{\gamma \in \gG (\B (\eta ))}f\left(\frac{1}{L}\gamma \right) \right| \notag \\ 
=&\frac{1}{L^{6g-6+2r}} \left| \sum_{\phi \in \gG/\Stab_\gG (\eta )}f\left( \frac{1}{L}\phi (\eta )\right) - s_\gG (\eta )\sum_{\phi \in \gG/\Stab_\gG (\B (\eta ))}f\left( \frac{1}{L}\phi (\B (\eta ))\right) \right| \notag \\ 
=&\frac{1}{L^{6g-6+2r}} \left| \sum_{\phi \in \gG/\Stab_\gG (\eta )}f\left( \frac{1}{L}\phi (\eta )\right) - \sum_{\phi \in \gG/\Stab_\gG (\eta )}f\left( \frac{1}{L}\B ( \phi (\eta ))\right) \right| \notag  \\ 
\leq &\frac{1}{L^{6g-6+2r}}\sum_{\phi \in \gG/\Stab_\gG (\eta )} \left| f\left( \frac{1}{L}\phi (\eta )\right) - f\left( \frac{1}{L}\B ( \phi (\eta ))\right) \right| .\tag{$\dagger$}
\end{align}

Let $\supp (f)$ denote the support of $f$.
Since $\ell_\SC \: \SC_K (\gS) \rightarrow \RRR$ is continuous, the compact set $\supp (f)$ is included in $\ell_\SC^{-1} ([0,D] ) $ for some $D>0$.

Note that for any $\mu\in \SC_K(\gS)$, $\mu\in \ell_\SC^{-1} ([0,D] ) $ as soon as $\B (\mu )\in \ell_\SC^{-1} ([0,D] ) $.
Hence, for $\phi \in \gG$, if $\frac{1}{L} \phi (\eta )\not \in \ell_\SC^{-1}([0,D])$, then $\frac{1}{L} \B (\phi (\eta )) \not \in \ell_\SC^{-1}([0,D])$, and so we have
\[ f\left( \frac{1}{L}\phi (\eta ) \right) =0= f\left( \frac{1}{L} \B (\phi (\eta ))\right) .\]
This concludes that in the last part ($\dagger$) of the above inequality, it is enough to consider the sum taken over $\phi$ belonging to
\[ \Phi_L =\left\{ \phi\in \gG/\Stab_\gG (\eta ) \; \middle| \; \frac{1}{L}\phi (\eta )\in \ell_\SC^{-1} ([0,D])\right\}.\] 

We observe that
\[ m_\eta^L (\ell_\SC^{-1}([0,D]  ))
=\frac{1}{L^{6g-6+2r}}\# \left\{ \eta' \in \gG (\eta )\ \middle| \ \ell_\SC \left( \frac{1}{L}\eta '\right)\leq D \right\}   = \frac{1}{L^{6g-6+2r}}\# \Phi_L \]
and by Theorem \ref{thm:counting subgroup use length},
\begin{align*}
m_\eta^L (\ell_\SC^{-1}([0,D]  ))&=D^{6g-6+2r}\frac{1}{(DL)^{6g-6+2r}}\# \{ \eta' \in \gG (\eta )\mid \ell_\SC ( \eta ')\leq DL \}\\
&\underset{L\rightarrow \infty }{\longrightarrow }D^{6g-6+2r}s_\gG (\eta )\mathfrak{c}^\Gamma_{g,r}(\B (\eta )) m_\Thu (\ell^{-1}([0,1])) <\infty .
\end{align*}
Hence, $m_\eta^L (\ell_\SC^{-1}([0,D]  ))$ is uniformly bounded above by some constant $M>0$.
We note that since any compact subset of $\SC_K (\gS)$ is included in $\ell_\SC^{-1}([0,D'])$ for some $D'>0$, $m_\eta^L$ is a locally finite measure on $\SC_K(\gS)$.

Take any $\varepsilon >0$. Take any metric function $d$ on $\SC_K (\gS)$ compatible with the topology. Since $f$ is uniformly continuous, we can take $\delta>0$ such that for any $x,y\in \SC_K(\gS)$ with $d(x,y)<\delta$, we have $|f(x)-f(y)|<\varepsilon/ M$.
Now, to apply to the inequality ($\dagger$), we want to see that
\[ \lim_{L\rightarrow \infty }\sup_{\phi \in \Phi_L } d\left( \frac{1}{L}\phi (\eta ), \frac{1}{L}\B (\phi (\eta ) )\right) =0.\]

Let us prove it by contradiction. To do so, suppose that there exist a positive constant $\tau>0$, a sequence $L_n\rightarrow \infty$ and $\phi_n\in \Phi_{L_n}$ such that for any $n\in \NN$ we have
\[ d\left( \frac{1}{L_n}\phi_n (\eta ), \frac{1}{L_n}\B (\phi_n (\eta ) )\right) >\tau.\]
Note that $\phi_n (\frac{1}{L_n} \eta )\in \ell_\SC^{-1}([0,D] )$ for all $n$. We aim to show that $\phi_n (\frac{1}{L_n} \eta )$ has a converging subsequence; however,$\ell_\SC^{-1}([0,D] )$ is not compact since $\ell_\SC (c\eta_G )=0$ for any $c>0$.

Now, we recall that we have the area functional $\Area \: \SC_K (\gS)\rightarrow \RRR$, which is continuous and $\RRR$-linear.
For any $\mu\in \SC_K (\gS)\setminus \GC_K (\gS)$, we have $\Area (\mu )>0$. Hence,
\[ \ell_A:=\ell_\SC +\Area \: \SC _K(\gS )\rightarrow \RRR \]
is a positive continuous $\RRR$-linear functional, which implies that $\ell_A^{-1}([0, T])$ is compact for any $T\geq 0$.
Since the area of each convex core is preserved by the action of $\Map (\gS)$, we have
\[ \ell_A \left( \phi_n \left( \frac{1}{L_n} \eta \right) \right)=\ell_\SC \left(\phi_n \left(\frac{1}{L_n} \eta \right)\right) +\Area \left(\phi_n \left(\frac{1}{L_n} \eta \right)\right)\leq D+\frac{1}{L_n} \Area (\eta )  .\]

Hence, the sequence $\phi_n (\frac{1}{L_n}\eta )$ is included in the compact set $\ell_A^{-1}([0,T])$ for some $T>0$, indicating that $\phi_n (\frac{1}{L_n} \eta )$ has a converging subsequence $\mu_n$. We denote the limit of $\mu_n$ by $\mu$.
Since 
\[ \Area \left(\phi_n \left(\frac{1}{L_n} \eta \right)\right)=\frac{1}{L_n} \Area (\eta ) \rightarrow 0 \quad (n\rightarrow \infty ),\]
we have $\Area (\mu )=0$, which implies that $\mu\in \GC_K (\gS)$.
Then, by the continuity of $\B$,
\[ 0< \tau \leq \lim_{n\rightarrow \infty }d(\mu_n , \B (\mu_n ) )=d(\mu , \B (\mu ))=d(\mu ,\mu )=0,\]
a contradiction.

From the above, for a sufficiently large $L$, we have
\[ d\left( \frac{1}{L}\phi (\eta ), \frac{1}{L}\B (\phi (\eta ) )\right) <\delta \]
for any $\phi \in \Phi_L$. Hence,
\begin{align*}
&\frac{1}{L^{6g-6+2r}}\sum_{\phi \in \gG/\Stab_\gG (\eta )} \left| f\left( \frac{1}{L}\phi (\eta )\right) - f\left( \frac{1}{L}\B ( \phi (\eta ))\right) \right| \\
=&\frac{1}{L^{6g-6+2r}}\sum_{\phi \in \Phi_L} \left| f\left( \frac{1}{L}\phi (\eta )\right) - f\left( \frac{1}{L}\B ( \phi (\eta ))\right) \right| \\
\leq & \frac{1}{L^{6g-6+2r}}\# \Phi_L \cdot \frac{\varepsilon}{M} \\
=& m_{\eta} ^L (\ell_\SC^{-1} ([0,D] )) \cdot \frac{\varepsilon}{M}\\
<& \varepsilon.
\end{align*}
This completes the proof.
\end{proof}

Similar to the proof of Corollary \ref{cor:ES}, we can establish the following corollary:

\begin{corollary}\label{cor:counting eta}
Under the conditions described above, for any positive homogeneous continuous function $F\: \SC_K(\gS) \rightarrow \RRR$, we have
\[ \lim_{L\rightarrow \infty}\frac{\# \{ \eta' \in \Gamma (\eta )\mid F (\eta' )\leq L\}}{L^{6g-6+2r}}=s_\gG (\eta ) \mathfrak{c}^\Gamma_{g,r}(\B (\eta )) m_\Thu (F^{-1}([0,1])).\]
\end{corollary}

\begin{remark}\label{rem:counting and area}
We must be cautious regarding the positivity of the function $F$ since $\ell_\SC (\eta _G)=0$, which implies that $\ell_\SC$ is not positive on $\SC_K(\gS)$.
Generally, when $F$ is equal to $F_0 \circ \B$ for any positive homogeneous continuous function $F_0\: \GC_K (\gS)\rightarrow \RRR$, $F$ is not a positive function on $\SC_K(\gS)$.

However, we can consider the function
\[ \Area + F_0\circ \B,\]
which is homogeneous, continuous, and positive on $\SC_K(\gS)$. In fact, for any non-zero $\mu \in \SC_K(\gS) $, $F_0\circ \B (\mu )>0$ if $\mu \in \GC_K(\gS)$, and $\Area (\mu )>0$ if $\mu\in \SC_K(\gS) \setminus \GC_K (\gS)$.

Let's consider the asymptotic formula for $\Area + F_0\circ \B$. For any $\eta'\in \gG (\eta)$, we have
\[  (\Area + F_0\circ \B ) (\eta ') =\Area (\eta ) +F_0\circ \B (\eta ') \]
since the action of $\Map (\gS)$ preserves the area of each convex core. The constant $\Area (\eta )$ does not influence the limit of the asymptotic formula, that is,
\begin{align*}
&\lim_{L\rightarrow \infty}\frac{\# \{ \eta' \in \Gamma (\eta )\mid (\Area + F_0\circ \B ) (\eta ')\leq L\}}{L^{6g-6+2r}} \\
=&\lim_{L\rightarrow \infty}\frac{\# \{ \eta' \in \Gamma (\eta )\mid F_0\circ \B (\eta ')\leq L-\Area (\eta )\}}{L^{6g-6+2r}}\\
=&\lim_{L\rightarrow \infty} \left( \frac{L-\Area (\eta )}{L}\right)^{6g-6+2r} \frac{\# \{ \eta' \in \Gamma (\eta )\mid F_0\circ \B (\eta ' )\leq L-\Area (\eta )\}}{(L-\Area (\eta ) )^{6g-6+2r}}\\
=&\lim_{L\rightarrow \infty}\frac{\# \{ \eta' \in \Gamma (\eta )\mid F_0\circ \B (\eta ' )\leq L\}}{L^{6g-6+2r}}.
\end{align*}
We also see that
\begin{align*}
m_\Thu ((\Area + F_0\circ \B )^{-1} ([0,1]) )
&= m_\Thu ((\Area + F_0\circ \B)^{-1} ([0,1]) \cap \ML (\gS ))\\
&=m_\Thu ( (F_0\circ \B )^{-1} ([0,1] )).
\end{align*}
Hence, we can formulate the following corollary.
\end{remark}

\begin{corollary}\label{cor:count general subset currents}
For any homogeneous continuous function $F\: \SC_K (\gS) \rightarrow \RRR$ that is positive on $\GC_K (\gS)$, we have
\[ \lim_{L\rightarrow \infty}\frac{\# \{ \eta' \in \Gamma (\eta )\mid F (\eta' )\leq L\}}{L^{6g-6+2r}}=s_\gG (\eta ) \mathfrak{c}^\Gamma_{g,r}(\B (\eta )) m_\Thu (F^{-1}([0,1])).\]
\end{corollary}

Given that $m_\Thu (F^{-1}([0,1]))$ is determined by the restriction of $F$ to $\ML (\gS)$, the constant $m_\Thu (F^{-1}([0,1]))$ can appear in the context of geodesic currents, as seen in Corollary \ref{cor:ES}.

This corollary leads to the asymptotic formula for weighted sum of conjugacy classes of $\Sub (G)$.
Let $J =a_1[H_1] +\cdots +a_m [H_m]$.
For a function $F\: \SC_K (\gS) \rightarrow \RRR$ and any $\phi (J) \in \gG (J )$ we define
\[ F(\phi (J)) =F(\phi (\eta ))=F(a_1\eta_{\phi (H_1)}+ \cdots + a_m \eta_{\phi (H_m)}) .\]
We can then formulate the following corollary, extending the first asymptotic formula in Theorem \ref{thm:counting subgroup use length}.

\begin{corollary}\label{cor:count general subgroups}
For any homogeneous continuous function $F\: \SC_K (\gS) \rightarrow \RRR$ that is positive on $\GC_K (\gS)$, the following holds:
\[ \lim_{L\rightarrow \infty}\frac{\# \{ J' \in \Gamma (J )\mid F (J' )\leq L\}}{L^{6g-6+2r}}=s_\gG (J ) \mathfrak{c}^\Gamma_{g,r}(\B (\eta ) ) m_\Thu (F^{-1}([0,1])).\]
The constant $s_\gG (J )$ originates from Lemma \ref{lem:B(H)/H is finite in general}.
\end{corollary}
\begin{proof}
Recall the relation
\[ \Stab_\gG (\B (\eta ))/ \Stab_\gG (J ) \cong \Stab_\gG (\B(\eta ))/ \Stab_\gG (\eta ) \times \Stab_\gG (\eta )/\Stab_\gG (J ) \]
and
\[ s_\gG (J)=s_\gG (\eta )\cdot  \# \Stab_\gG (\eta )/\Stab_\gG (J ).\]
Hence, we have
\begin{align*}
&\frac{1}{L^{6g-6+2r}}\# \{ J' \in \Gamma (J )\mid F (J' )\leq L\} \\
=&\frac{1}{L^{6g-6+2r}}\# \{ (\phi, \psi )\in \gG /\Stab_\gG (\eta )\times \Stab_\gG (\eta )/\Stab_\gG (J ) \mid F (\phi \circ \psi (J) )\leq L\} \\
=&\frac{1}{L^{6g-6+2r}}\# \{ (\phi, \psi )\in \gG /\Stab_\gG (\eta )\times \Stab_\gG (\eta )/\Stab_\gG (J ) \mid F (\phi \circ \psi (\eta ) )\leq L\} \\
=&\sum_{\psi \in \Stab_\gG (\eta )/\Stab_\gG (J )} \frac{1}{L^{6g-6+2r}}\# \{ \phi \in \gG/\Stab_\gG (\eta ) \mid F (\phi (\eta ) )\leq L \} \\
=&\# \Stab_\gG (\eta )/\Stab_\gG (J )\cdot \frac{1}{L^{6g-6+2r}}\# \{ \eta' \in \gG (\eta ) \mid F (\eta ' )\leq L \}\\
\underset{L\rightarrow \infty }{\longrightarrow } &\# \Stab_\gG (\eta )/\Stab_\gG (J ) \cdot s_\gG (\eta ) \mathfrak{c}^\Gamma_{g,r}(\B (\eta )) m_\Thu (F^{-1}([0,1])) \\
=&s_\gG (J ) \mathfrak{c}^\Gamma_{g,r}(\B (\eta )) m_\Thu (F^{-1}([0,1])) .
\end{align*}
This completes the proof.
\end{proof}

\begin{example}[Continuous functional on $\SC (\gS)$]
Let's assume that $\gS$ is a closed hyperbolic surface. In this case, we do not have to consider the subspace $\SC_K (\gS)$ of $\SC (\gS)$ for some compact $K\subset \gS$.
We aim to identify a function $F\: \SC (\gS)\rightarrow \RRR$ satisfying the conditions in Corollary \ref{cor:count general subset currents} and \ref{cor:count general subgroups}, which cannot be expressed as $F_0 \circ \B$ for any homogeneous continuous function $F_0\: \GC (\gS) \rightarrow \RRR$.
Remark that if $F=F_0 \circ \B$, then for any $\mu \in \GC (\gS)$ we have
\[ F(\mu )=(F_0\circ \B)(\mu)=F_0 (\mu),\]
which implies that $F_0$ coincides with the restriction $F|_{\GC (\gS)}$ of $F$ to $\GC (\gS)$.

We also note that even if $F=F_0\circ \B$, $\Area+ F$ cannot be expressed as $F_1\circ \B$ for any homogeneous continuous function $F_1\: \GC (\gS) \rightarrow \RRR$; however, the area functional $\Area$ does not influence the counting formula (see Remark \ref{rem:counting and area}). Hence, we seek another example.

Previous work \cite{Sas22} has introduced several continuous functionals on $\SC (\gS)$.
We focus on the generalized intersection number functional $i_\SC \: \SC (\gS) \times \SC (\gS )\rightarrow \RRR$, which is introduced in \cite[Theorem 5.39]{Sas22}.
We are going to construct a functional from $\SC (\gS)$ to $\RRR$ with the desired property using $i_\SC$ (see Proposition \ref{prop:i_SCnot=}).
For the reader's convenience, we briefly review the definition and properties of $i_\SC$.

For $H,K \in \Sub (G )$ we consider the diagonal action of $G$ on $G/H \times G/K$ and the associated quotient set $G\backslash (G/H\times G/K)$. The intersection number $i(C_H,C_K)$ between $C_H$ and $C_K$ is defined as the number of the equivalence classes, $[g_1H, g_2K]\in G\backslash (G/H\times G/K)$, that satisfy the condition that $g_1\CH (\gL (H) )\cap g_2 \CH (\gL (K))$ forms a non-empty compact set.
When $H$ and $K$ are cyclic, $i(C_H, C_K)$ coincides with the (geometric) intersection number of two closed geodesics $C_H$ and $C_K$ on $\gS$.

We note that $i(C_H,C_K)$ is equal to the number of contractible components of the fiber product 
\[ C_H\times _\gS C_K:= \{ (x,y)\in C_H\times C_K \mid p_H (x)=p_K (y) \} \]
with respect to the canonical projections $p_H\: C_H\rightarrow \gS$ and $p_K\: C_K\rightarrow \gS$.

The intersection number functional $i_\SC$ is a continuous, symmetric, $\Map (\gS)$-invariant, and $\RRR$-bilinear functional satisfying that for any $H,K\in \Sub (G)$ we have
\[ i_\SC (\eta _H,\eta_K )=i (C_H, C_K ).\]
The restriction of $i_\SC$ to $\GC (\gS )\times \GC (\gS)$, denoted by $i_\GC$, coincides with the continuous extension of the intersection number of closed geodesics introduced by Bonahon \cite{Bon86}.
When we fix $\mu \in \SC (\gS )$, we can obtain the continuous $\RRR$-linear functional
\[ i_\mu =i_\SC (\mu, \cdot )\: \SC (\gS ) \rightarrow \RRR .\]
In Proposition \ref{prop:i_SCnot=}, we will prove that there exists $\mu \in \SC (\gS)\setminus \GC (\gS)$ such that $i_\mu \not= i_\mu|_{\GC(\gS)}\circ \B$.

We review the argument in \cite[Example 5.15, Theorem 7.4]{Sas22}. Consider a cyclic subgroup $H\in \Sub (G)$ and any non-cyclic $K\in \Sub (G)$.
If $g_1\CH (\gL (H ) ) \cap g_2\CH (\gL (K))$ is a non-empty compact set, then this intersection is a geodesic segment.
The endpoints of this segment arise from the intersection points between $C_H$ and $\partial C_K$ in $\gS$. Consequently, we establish the equality:
\[ i(C_H, C_K )=\frac{1}{2}\sum _{c\in \partial C_K }i (C_H , c) =i (\eta_H , \B (\eta_K )).\]
This implies that if $H$ is cyclic, then we have
\[ i_{\eta_H } = i_\GC (\eta_H, \B (\cdot )),\]
that is, $i_{\eta_H}\: \SC (\gS )\rightarrow \RRR$ equals the composition of $\B$ and $i_\GC (\eta_H, \cdot ) \: \GC (\gS )\rightarrow \RRR$, which does not meet our objective.

Recall that $\mu \in \GC (\gS )$ is said to be \emph{filling} if $i_\GC (\mu ,\nu ) >0$ for every non-zero $\nu \in \GC (\gS)$.
We focus on $\mu \in \SC (\gS)$ with the property that $\B (\mu ) $ is filling. From the above argument, for any $\nu \in \GC (\gS)$ we have
\[ i_\SC (\mu ,\nu )=i_\GC (\B (\mu ) , \nu ) >0,\]
implying that the restriction $i_\mu |_{\GC (\gS)} \: \GC (\gS) \rightarrow \RRR$ is a positive continuous $\RRR$-linear functional.
Consequently, $i_\mu$ satisfies the conditions in Corollary \ref{cor:count general subgroups}.

We emphasize the following inequality: for any $\mu,\nu \in \SC (\gS)$, we have
\[ i_\SC (\mu ,\nu )\leq i_\GC (\B (\mu ),\B (\nu ) ).\]
This inequality is justified because any non-empty compact intersection $\CH (\gL (H))\cap \CH (\gL (K))$ forms a polygon with at least four vertices for non-cyclic $H,K\in \Sub (G)$ (see \cite[Theorem 7.4]{Sas22} for detail).
However, we lack a lower bound for $i_\SC (\mu, \nu )$.
It is worth mentioning that the intersection number $i_\SC (\mu ,\nu )$ can be zero, even when $i_\GC (\B(\mu ) ,\B (\nu ))>0$.

\begin{figure}[h]
\centering
\includegraphics{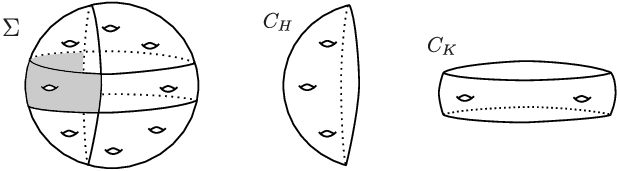}
\caption{In the left of the figure, $p_H(C_H)$ and $p_K(C_K)$ are described as subsurfaces of the closed surface $\gS$ of genus $8$. Their intersection $p_H(C_H)\cap p_K (C_K)$, which is blacked out, is a surface of genus $1$ with one boundary component.}\label{fig:intersection_number_of_simple_compact_surfaces}
\end{figure}

Actually, in the case of Figure \ref{fig:intersection_number_of_simple_compact_surfaces}, the intersection $p_H(C_H)\cap p_K (C_K)$, which is homeomorphic to the fiber product $C_H\times _\gS C_K$, is not contractible. Hence, $i_\SC (\eta_H,\eta_K)=0$. Nevertheless, $i_\GC (\B (\eta_H),\B (\eta_K))=\frac{1}{2}\cdot \frac{1}{2}\cdot 4=1>0$.
By using this example, we can obtain the following proposition:
\end{example}

\begin{proposition}\label{prop:i_SCnot=}
Let $\gS$ be a closed hyperbolic surface. Then, there exists $\mu \in \SC (\gS)\setminus \GC (\gS)$ such that $\B (\mu)$ is filling and 
\[
i_\mu \not= i_\mu|_{\GC(\gS)}\circ \B.
\]
As a result, the $\RRR$-linear functional $i_\mu$, which is positive on $\GC (\gS)$, cannot be expressed as $F_0 \circ \B$ for any homogeneous continuous function $F_0\: \GC (\gS) \rightarrow \RRR$.
\end{proposition}
\begin{proof}
Take $\mu_0 \in \SC (\gS)$ such that $\B (\mu_0)$ is filling. Note that $\mu_0$ may possibly belong to $\GC (\gS)$.
Using the subgroups $H$ and $K$ of $G$ from the above example, define $\mu:=\mu_0+ \eta_H$, which belongs to $\SC (\gS )\setminus \GC (\gS)$.
Then we have
\begin{align*}
  i_\mu (\eta_K)=&i_\SC (\mu +\eta_H, \eta_K)\\
  =&i_\SC (\mu,\eta_K)+i_\SC (\eta_H,\eta_K)\\
  =&i_\SC (\mu,\eta_K)\\
  \leq &i_\GC (\B(\mu),\B(\eta_K))\\
  =&i_\SC (\mu ,\B (\eta_K))\\
  <&i_\SC (\mu, \B (\eta_K))+i_\GC (\B(\eta_H), \B (\eta_K))\\
  =&i_\SC (\mu, \B (\eta_K))+i_\SC (\eta_H,\B (\eta_K))\\
  =&i_\mu (\B (\eta_K)).
\end{align*}
This implies that $i_\mu \not= i_\mu |_{\GC (\gS)}\circ \B$.
\end{proof}

For a subset current $\mu \in \SC (\gS)$ satisfying the condition of the above proposition, we obtain the asymptotic formula:
\[ \lim_{L\rightarrow \infty}\frac{\# \{ \eta' \in \Gamma (\eta )\mid i_\mu (\eta' )\leq L\}}{L^{6g-6}}=s_\gG (\eta ) \mathfrak{c}^\Gamma_{g,r}(\B (\eta ) ) m_\Thu (i_\mu ^{-1}([0,1])),\]
which cannot be deduced from Theorem \ref{thm:counting subgroup use length}.

From the above, we can formulate the following theorem.
\begin{theorem}
Let $\gS$ be a closed hyperbolic surface of genus $g\geq 2$.
For any $\mu \in \SC (\gS)$, if $\B (\mu )$ is filling, then we have
\[ \lim_{L\rightarrow \infty}\frac{\# \{ \eta' \in \Gamma (\eta )\mid i_\mu (\eta' )\leq L\}}{L^{6g-6}}=s_\gG (\eta ) \mathfrak{c}^\Gamma_{g,r}(\B (\eta ) ) m_\Thu (i_\mu ^{-1}([0,1])),\]
where $i_\mu= i_\SC (\mu , \cdot )$.
Specifically, we have
\[ \sup_{\phi \in \Gamma } i_\SC (\mu ,\phi (\eta ) ) =\infty.\]
\end{theorem}

Finally, we present an example of $\mu \in \SC (\gS) \setminus \GC (\gS)$ satisfying the condition of $\B (\mu )$ being filling.
Consider a filling geodesic current $\nu$ represented as
\[ \nu =\eta_{c_1}+\cdots + \eta_{c_k} ,\]
where $c_1,\dots c_k$ are closed geodesics on $\gS$.
By \cite[Theorem 7.9]{Sas22} and the preceding argument of the theorem, we can find a non-cyclic $H_i\in \Sub (G)$ such that $\B (\eta_{H_i} )=\eta_{c_i}$ for $i=1,\dots ,k$.
Thus, $\B (\eta_{H_1}+\cdots +\eta_{H_k}) =\nu$.

We note that if $c_i$ is simple, then $H_i$ can be obtained by cutting $\gS$ along $c_i$.
By applying this procedure, we can construct $\mu \in \SC (\gS)$ such that $\B (\mu )$ is filling, expressed as
\[ \mu =\eta_{H_1}+\cdots +\eta_{H_j} ,\]
where $H_1,\dots ,H_j\in \Sub (G)$ are non-cyclic.

\end{document}